\documentclass{amsart}
\usepackage{amssymb}
\usepackage{latexsym}
\usepackage{hyperref}

\date{\today}

\newcommand{\Z}{{\mathbb Z}}
\newcommand{\R}{{\mathbb R}}
\newcommand{\C}{{\mathbb C}}

\newcommand{\Q}{{\mathbb Q}}
\newcommand{\D}{{\mathbb D}}

\newtheorem{theorem}{Theorem}[section]

\newtheorem{lemma}[theorem]{Lemma}
\newtheorem{prop}[theorem]{Proposition}
\newtheorem{coro}[theorem]{Corollary}
\newtheorem{defi}[theorem]{Definition}
\begin{document}
\title[Strictly Ergodic Subshifts and Associated Operators]{Strictly Ergodic Subshifts and\\Associated Operators}

\author{David Damanik}

\address{Mathematics 253--37, California Institute of Technology, Pasadena, CA 91125,
USA}

\email{\href{mailto:damanik@caltech.edu}{damanik@caltech.edu}}

\urladdr{\href{http://www.math.caltech.edu/people/damanik.html}{http://www.math.caltech.edu/people/damanik.html}}

\thanks{This work was supported in part by NSF grant DMS--0500910}

\begin{abstract}
We consider ergodic families of Schr\"odinger operators over base dynamics given by
strictly ergodic subshifts on finite alphabets. It is expected that the majority of these
operators have purely singular continuous spectrum supported on a Cantor set of zero
Lebesgue measure. These properties have indeed been established for large classes of
operators of this type over the course of the last twenty years. We review the mechanisms
leading to these results and briefly discuss analogues for CMV matrices.
\end{abstract}

\dedicatory{Dedicated to Barry Simon on the occasion of his 60th birthday.}

\maketitle

\section{Introduction}\label{Sec1}

When I was a student in the mid-1990's at the Johann Wolfgang Goethe Universit\"at in
Frankfurt, my advisor Joachim Weidmann and his students and postdocs would meet in his
office for coffee every day and discuss mathematics and life. One day we walked in and
found a stack of preprints on the coffee table. What now must seem like an ancient
practice was not entirely uncommon in those days: In addition to posting preprints on the
archives, people would actually send out hardcopies of them to their peers around the
world.

In this particular instance, Barry Simon had sent a series of preprints, all dealing with
singular continuous spectrum. At the time I did not know Barry personally but was well
aware of his reputation and immense research output. I was intrigued by these preprints.
After all, we had learned from various sources (including the Reed-Simon books!) that
singular continuous spectrum is sort of a nuisance and something whose absence should be
proven in as many cases as possible. Now we were told that singular continuous spectrum
is generic?

Soon after reading through the preprint series it became clear to me that my thesis topic
should have something to do with this beast: singular continuous spectrum.
Coincidentally, only a short while later I came across a beautifully written paper by
S\"ut\H{o} \cite{s5} that raised my interest in the Fibonacci operator. I had studied
papers on the almost Mathieu operator earlier. For that operator, singular continuous
spectrum does occur, but only in very special cases, that is, for special choices of the
coupling constant, the frequency, or the phase. In the Fibonacci case, however, singular
continuous spectrum seemed to be the rule. At least there was no sensitive dependence on
the coupling constant or the frequency as I learned from \cite{bist,s5,s6}.

Another feature, which occurs in the almost Mathieu case, but only at special coupling,
seemed to be the rule for the Fibonacci operator: zero-measure spectrum.

So I set out to understand what about the Fibonacci operator was responsible for this
persistent occurrence of zero-measure singular continuous spectrum. Now, some ten years
later, I still do not really understand it. In fact, as is always the case, the more you
understand (or think you understand), the more you realize how much else is out there,
still waiting to be understood.\footnote{Most recently, I have come to realize that I do
not understand why the Lyapunov exponent vanishes on the spectrum, even at large
coupling. Who knows what will be next...}

Thus, this survey is meant as a snapshot of the current level of understanding of things
related to the Fibonacci operator and also as a thank-you to Barry for having had the
time and interest to devote a section or two of his OPUC book to subshifts and the
Fibonacci CMV matrix. Happy Birthday, Barry, and thank you for being an inspiration to so
many generations of mathematical physicists!

\section{Strictly Ergodic Subshifts}\label{Sec2}

In this section we define strictly ergodic subshifts over a finite alphabet and discuss
several examples that have been studied from many different perspectives in a great
number of papers.

\subsection{Basic Definitions}

We begin with the definitions of the the basic objects:

\begin{defi}[full shift]
Let $\mathcal{A}$ be a finite set, called the alphabet. The two-sided infinite sequences
with values in $\mathcal{A}$ form the full shift $\mathcal{A}^\Z$. We endow $\mathcal{A}$
with the discrete topology and the full shift with the product topology.
\end{defi}

\begin{defi}[shift transformation]
The shift transformation $T$ acts on the full shift by $[T \omega]_n = \omega_{n+1}$.
\end{defi}

\begin{defi}[subshift]
A subset $\Omega$ of the full shift is called a subshift if it is closed and
$T$-invariant.
\end{defi}

Thus, our base dynamical systems will be given by $(\Omega,T)$, where $\Omega$ is a
subshift and $T$ is the shift transformation. This is a special class of topological
dynamical systems that is interesting in its own right. Basic questions regarding them
concern the structure of orbits and invariant (probability) measures. The situation is
particularly simple when orbit closures and invariant measures are unique:

\begin{defi}[minimality]
Let $\Omega$ be a subshift and $\omega \in \Omega$. The orbit of $\omega$ is given by
$\mathcal{O}_\omega = \{ T^n \omega : n \in \Z \}$. If $\mathcal{O}_\omega$ is dense in
$\Omega$ for every $\omega \in \Omega$, then $\Omega$ is called minimal.
\end{defi}

\begin{defi}[unique ergodicity]
Let $\Omega$ be a subshift. A Borel measure $\mu$ on $\Omega$ is called $T$-invariant if
$\mu(T(A)) = \mu(A)$ for every Borel set $A \subseteq \Omega$. $\Omega$ is called
uniquely ergodic if there is a unique $T$-invariant Borel probability measure on
$\Omega$.
\end{defi}

By compactness of $\Omega$, the set of $T$-invariant Borel probability measure on
$\Omega$ is non-empty. It is also convex and the extreme points are exactly the ergodic
measures, that is, probability measures for which $T(A) = A$ implies that either $\mu(A)
= 0$ or $\mu(A) = 1$. Thus, a subshift is uniquely ergodic precisely when there is a
unique ergodic measure on it.

We will focus our main attention on subshifts having both of these properties. For
convenience, one often combines these two notions into one:

\begin{defi}[strict ergodicity]
A subshift $\Omega$ is called strictly ergodic if it is both minimal and uniquely
ergodic.
\end{defi}

\subsection{Examples of Strictly Ergodic Subshifts}

Let us list some classes of strictly ergodic subshifts that have been studied by a
variety of authors and from many different perspectives (e.g., symbolic dynamics, number
theory, spectral theory, operator algebras, etc.) in the past.

\subsubsection{Subshifts Generated by Sequences}

Here we discuss a convenient way of defining a subshift starting from a sequence $s \in
\mathcal{A}^\Z$. Since $\mathcal{A}^\Z$ is compact, $\mathcal{O}_s$ has a non-empty set
of accumulation points, denoted by $\Omega_s$. It is readily seen that $\Omega_s$ is
closed and $T$-invariant. Thus, we call $\Omega_s$ the \textit{subshift generated by
$s$}. Naturally, we seek conditions on $s$ that imply that $\Omega_s$ is minimal or
uniquely ergodic.

Every word $w$ (also called block or string) of the form $w = s_m \ldots s_{m+n-1}$ with
$m \in \Z$ and $n \in \Z^+ = \{ 1,2,3,\ldots \}$ is called a \textit{subword of} $s$ (of
\textit{length} $n$, denoted by $|w|$). We denote the set of all subwords of $s$ of
length $n$ by $\mathcal{W}_s(n)$ and let
$$
\mathcal{W}_s = \bigcup_{n \ge 1} \mathcal{W}_s(n).
$$
If $w \in \mathcal{W}_s(n)$, let $\cdots < m_{-1} < 0 \le  m_0 < m_1 < \cdots$ be the
integers $m$ for which $s_m \ldots s_{m+n-1} = w$. The sequence $s$ is called
\textit{recurrent} if $m_n \to \pm \infty$ as $n \to \pm \infty$ for every $w \in
\mathcal{W}_s$. A recurrent sequence $s$ is called \textit{uniformly recurrent} if
$(m_{j+1} - m_j)_{j \in \Z}$ is bounded for every $w \in \mathcal{W}_s$. Finally, a
uniformly recurrent sequence $s$ is called \textit{linearly recurrent} if there is a
constant $C < \infty$ such that for every $w \in \mathcal{W}_s$, the gaps $m_{j+1} - m_j$
are bounded by $C |w|$.

We say that $w \in \mathcal{W}_s$ \textit{occurs in $s$ with a uniform frequency} if
there is $\mathrm{d}_s(w) \ge 0$ such that, for every $k \in \Z$,
$$
\mathrm{d}_s(w) = \lim_{n \to \infty} \frac{1}{n} \left| \{ m_j \}_{j \in \Z} \cap
[k,k+n) \right|,
$$
and the convergence is uniform in $k$.

For results concerning the minimality and unique ergodicity of the subshift $\Omega_s$
generated by a sequence $s$, we recommend the book by Queff\'elec \cite{q}; see in
particular Section~IV.2. Let us recall the main findings.

\begin{prop}
If $s$ is uniformly recurrent, then $\Omega_s$ is minimal. Conversely, if $\Omega$ is
minimal, then every $\omega \in \Omega$ is uniformly recurrent. Moreover,
$\mathcal{W}_{\omega_1} = \mathcal{W}_{\omega_2}$ for every $\omega_1,\omega_2 \in
\Omega$.
\end{prop}

The last statement permits us to define a set $\mathcal{W}_\Omega$ for any minimal
subshift $\Omega$ so that $\mathcal{W}_\Omega = \mathcal{W}_\omega$ for every $\omega \in
\Omega$.

\begin{prop}
Let $s$ be recurrent. Then, $\Omega_s$ is uniquely ergodic if and only if each subword of
$s$ occurs with a uniform frequency.
\end{prop}

As a consequence, $\Omega_s$ is strictly ergodic if and only if each subword $w$ of $s$
occurs with a uniform frequency $\mathrm{d}_s(w) > 0$.

An interesting class of strictly ergodic subshifts is given by those subshifts that are
generated by linearly recurrent sequences \cite{dur,lagple}:

\begin{prop}\label{lrse}
If $s$ is linearly recurrent, then $\Omega_s$ is strictly ergodic.
\end{prop}

\subsubsection{Sturmian Sequences}

Suppose $s$ is a uniformly recurrent sequence. We saw above that $\Omega_s$ is a minimal
subshift and all elements of $\Omega_s$ have the same set of subwords,
$\mathcal{W}_{\Omega_s} = \mathcal{W}_s$. Let us denote the cardinality of
$\mathcal{W}_s(n)$ by $p_s(n)$. The map $\Z^+ \to \Z^+, \; n \mapsto p_s(n)$ is called
the \textit{complexity function of} $s$ (also called \textit{factor, block} or
\textit{subword complexity function}).

It is clear that a periodic sequence gives rise to a bounded complexity function. It is
less straightforward that every non-periodic sequence gives rise to a complexity function
that grows at least linearly. This fact is a consequence of the following celebrated
theorem due to Hedlund and Morse \cite{mh1}:

\begin{theorem}
If $s$ is recurrent, then the following statements are equivalent:
\begin{itemize}

\item[{\rm (i)}] $s$ is periodic, that is, there exists $k$ such that $s_m = s_{m+k}$ for
every $m \in \Z$.

\item[{\rm (ii)}] $p_s$ is bounded, that is, there exists $p$ such that $p_s(n) \le p$
for every $n \in \Z^+$.

\item[{\rm (iii)}] There exists $n_0 \in \Z^+$ such that $p_s(n_0) \le n_0$.

\end{itemize}
\end{theorem}

\begin{proof}
The implications (i) $\Rightarrow$ (ii) and (ii) $\Rightarrow$ (iii) are obvious, so we
only need to show (iii) $\Rightarrow$ (i).

Let $R_s(n)$ be the directed graph with $p_s(n)$ vertices and $p_s(n+1)$ edges which is
defined as follows. Every subword $w \in \mathcal{W}_s(n)$ corresponds to a vertex of
$R_s(n)$. Every $\tilde w \in \mathcal{W}_s(n+1)$ generates an edge of $R_s(n)$ as
follows. Write $\tilde w = axb$, where $a,b \in \mathcal{A}$ and $x$ is a (possibly
empty) string. Then draw an edge from the vertex $ax$ to the vertex $xb$.

We may assume that $\mathcal{A}$ has cardinality at least two since otherwise the theorem
is trivial. Thus, $p_s(1) \ge 2 > 1$. Obviously, $p_s$ is non-decreasing. Thus, by
assumption (iii), there must be $1 \le n_1 < n_0$ such that $p_s(n_1) = p_s(n_1 + 1)$.
Consider the graph $R_s(n_1)$. Since $s$ is recurrent, there must be a directed path from
$w_1$ to $w_2$ for every pair $w_1,w_2 \in \mathcal{W}_s(n_1)$. On the other hand,
$R_s(n_1)$ has the same number of vertices and edges. It follows that $R_s(n_1)$ is a
simple cycle and hence $s$ is periodic of period $p_s(n_1)$.
\end{proof}

The graph $R_s(n)$ introduced in the proof above is called the \textit{Rauzy graph}
associated with $s$ and $n$. It is an important tool for studying (so-called)
combinatorics on words. This short proof of the Hedlund-Morse Theorem is just one of its
many applications.

\begin{coro}
If $s$ is recurrent and not periodic, then $p_s(n) \ge n+1$ for every $n \in \Z^+$.
\end{coro}

This raises the question whether aperiodic sequences of minimal complexity exist.

\begin{defi}
A sequence $s$ is called Sturmian if it is recurrent and satisfies $p_s(n) = n+1$ for
every $n \in \Z^+$.
\end{defi}

\noindent\textbf{Remarks.} (a) There are non-recurrent sequences $s$ with complexity
$p_s(n) = n+1$. For example, $s_n = \delta_{n,0}$. The subshifts generated by such
sequences are trivial and we therefore restrict our attention to recurrent sequences.\\
(b) We have seen that growth strictly between bounded and linear is impossible for a
complexity function. It is an interesting open problem to characterize the increasing
functions from $\Z^+$ to $\Z^+$ that arise as complexity functions.

\medskip

Note that a Sturmian sequence is necessarily defined on a two-symbol alphabet
$\mathcal{A}$. Without loss of generality, we restrict our attention to $\mathcal{A} = \{
0,1 \}$. The following result gives an explicit characterization of all Sturmian
sequences with respect to this normalization.

\begin{theorem}\label{sturmdesc}
A sequence $s \in \{ 0,1 \}^\Z$ is Sturmian if and only if there are $\theta \in (0,1)$
irrational and $\phi \in [0,1)$ such that either
\begin{equation}\label{sturmqp}
s_n = \chi_{[1-\theta,1)} (n \theta + \phi) \; \text{ or } \;  s_n = \chi_{(1-\theta,1]}
(n \theta + \phi)
\end{equation}
for all $n \in \Z$.
\end{theorem}

\noindent\textbf{Remark.} In \eqref{sturmqp}, we consider the $1$-periodic extension of
the function $\chi_{[1-\theta,1)} (\cdot)$ (resp., $\chi_{(1-\theta,1]} (\cdot)$) on
$[0,1)$.

\medskip

Each sequence of the form \eqref{sturmqp} generates a subshift. The following theorem
shows that the resulting subshift only depends on $\theta$.

\begin{theorem}
Assume $\theta \in (0,1)$ is irrational, $\phi \in [0,1)$, and $s_n = \chi_{[1-\theta,1)}
(n \theta + \phi)$. Then the subshift generated by $s$ is given by
$$
\Omega_s = \left\{ n \mapsto \chi_{[1-\theta,1)} (n \theta + \tilde \phi) : \tilde \phi
\in [0,1) \right\} \cup \left\{ n \mapsto \chi_{(1-\theta,1]} (n \theta + \tilde \phi) :
\tilde \phi \in [0,1) \right\}.
$$
Moreover, $\Omega_s$ is strictly ergodic.
\end{theorem}

Let us call a subshift \textit{Sturmian} if it is generated by a Sturmian sequence. We
see from the previous theorem that there is a one-to-one correspondence between
irrational numbers $\theta$ and Sturmian subshifts. We call $\theta$ the \textit{slope}
of the subshift.

\medskip

\noindent\textbf{Example} (Fibonacci case). The Sturmian subshift corresponding to the
inverse of the golden mean,
$$
\theta = \frac{\sqrt{5} - 1}{2},
$$
is called the \textit{Fibonacci subshift} and its elements are called \textit{Fibonacci
sequences}.

\medskip

An important property of Sturmian sequences is their hierarchical, or $S$-\textit{adic},
structure. That is, there is a natural level of hierarchies such that on each level,
there is a unique decomposition of the sequence into blocks of two types. The starting
level is just the decomposition into individual symbols. Then, one may pass from one
level to the next by a set of rules that is determined by the coefficients in the
continued fraction expansion of the slope $\theta$.

Let
\begin{equation}\label{cfexp}
\theta = \cfrac{1}{a_1 + \cfrac{1}{a_2 + \cfrac{1}{a_3 + \cdots}}}
\end{equation}
be the continued fraction expansion of $\theta$ with uniquely determined $a_k \in \Z^+$.
Truncation of this expansion after $k$ steps yields rational numbers $p_k/q_k$ that obey
\begin{alignat}{3}
\label{pkdef}
p_0 &= 0, &\quad    p_1 &= 1,   &\quad  p_k &= a_k p_{k-1} + p_{k-2},\\
\label{qkdef} q_0 &= 1, &     q_1 &= a_1, &       q_k &= a_k q_{k-1} + q_{k-2}.
\end{alignat}
These rational numbers are known to be best approximants to $\theta$. See Khinchin
\cite{khin} for background on continued fraction expansions.

We define words $(w_k)_{k \in \Z^+_0}$ over the alphabet $\{0,1\}$ as follows:
\begin{equation}\label{wkrec}
w_0 = 0, \quad w_1 = 0^{a_1 - 1} 1, \quad w_{k+1} = w_k^{a_{k+1}} w_{k-1} \; \mbox{ for }
k \ge 1.
\end{equation}

\begin{theorem}\label{partition}
Let $\Omega$ be a Sturmian subshift with slope $\theta$ and let the words $w_k$ be
defined by \eqref{cfexp} and \eqref{wkrec}. Then, for every $k \in \Z^+$, each $\omega
\in \Omega$ has a unique partition, called the $k$-partition of $\omega$, into blocks of
the form $w_k$ or $w_{k-1}$. In this partition, blocks of type $w_k$ occur with
multiplicity $a_{k+1}$ or $a_{k+1} + 1$ and blocks of type $w_{k-1}$ occur with
multiplicity one.
\end{theorem}

\begin{proof}[Sketch of proof.]
The first step is to use the fact that $p_k/q_k$ are best approximants to show that the
restriction of $\chi_{[1-\theta,1)} (n \theta)$ to the interval $[1,q_k]$ is given by
$w_k$, $k \in \Z^+$; compare \cite{bist}. The recursion \eqref{wkrec} therefore yields a
$k$-partition of $\chi_{[1-\theta,1)} (n \theta)$ on $[1,\infty)$. Since every $\omega
\in \Omega$ may be obtained as an accumulation point of shifts of this sequence, it can
then be shown that a unique partition of $\omega$ is induced; see \cite{dl1}. The
remaining claims follow quickly from the recursion \eqref{wkrec}
\end{proof}

\noindent\textbf{Example} (Fibonacci case, continued). In the Fibonacci case, $a_k = 1$
for every $k$. Thus, both $(p_k)$ and $(q_k)$ are sequences of Fibonacci numbers (i.e.,
$p_{k+1} = q_k = F_k$, where $F_0 = F_1 = 1$ and $F_{k+1} = F_k + F_{k-1}$ for $k \ge 1$)
and the words $w_k$ are obtained by the simple rule
\begin{equation}\label{fibwkrec}
w_0 = 0, \quad w_1 = 1, \quad w_k = w_{k-1} w_{k-2} \; \mbox{ for } k \ge 2.
\end{equation}
Thus, the sequence $(w_k)_{k \in \Z^+}$ is given by $1, \; 10, \; 101, \; 10110, \;
10110101, \ldots$, which may also be obtained by iterating the rule
\begin{equation}\label{fibsubst}
1 \mapsto 10 , \quad 0 \mapsto 1,
\end{equation}
starting with the symbol $1$.

\medskip

For the proofs omitted in this subsection and much more information on Sturmian sequences
and subshifts, we refer the reader to \cite{berstel,dl13,loth2,mh2}.

\subsubsection{Codings of Rotations}

Theorem~\ref{sturmdesc} shows that Sturmian sequences are obtained by coding an
irrational rotation of the torus according to a partition of the circle into two
half-open intervals. It is natural to generalize this and consider codings of rotations
with respect to a more general partition of the circle. Thus, let $[0,1) = I_1 \cup
\ldots \cup I_l$ be a partition into $l$ half-open intervals. Choosing numbers
$\lambda_1, \ldots, \lambda_l$, we consider the sequences
\begin{equation}\label{codofrot}
s_n = \sum_{j = 1}^l \lambda_j \chi_{I_j}(n \theta + \phi).
\end{equation}
Subshifts generated by sequences of this form will be said to be associated with codings
of rotations.

\begin{theorem}
Let $\theta \in (0,1)$ be irrational and $\phi \in [0,1)$. If $s$ is of the form
\eqref{codofrot}, then $\Omega_s$ is strictly ergodic. Moreover, the complexity function
satisfies $p_s(n) = a n + b$ for every $n \ge n_0$ and suitable integers $a,b,n_0$.
\end{theorem}

See \cite{h} for a proof of strict ergodicity and \cite{ab2} for a proof of the
complexity statement. In fact, the integers $a,b,n_0$ can be described explicitly; see
\cite[Theorem~10]{ab2}.

\subsubsection{Arnoux-Rauzy and Episturmian Subshifts}

Let us consider a minimal subshift $\Omega$ over the alphabet $\mathcal{A}_m = \{
0,1,2,\ldots,m-1\}$, where $m \ge 2$. A word $w \in \mathcal{W}_\Omega$ is called
\textit{right-special} (resp., \textit{left-special}) if there are distinct symbols $a,b
\in \mathcal{A}_m$ such that $wa,wb \in \mathcal{W}_\Omega$ (resp., $aw,bw \in
\mathcal{W}_\Omega$). A word that is both right-special and left-special is called
\textit{bispecial}. Thus, a word is right-special (resp., left-special) if and only if
the corresponding vertex in the Rauzy graph has out-degree (resp., in-degree) $\ge 2$.

Note that the complexity function of a Sturmian subshift obeys $p(n+1) - p(n) = 1$ for
every $n$, and hence for every length, there is a unique right-special factor and a
unique left-special factor, each having exactly two extensions.

Arnoux-Rauzy subshifts and episturmian subshifts relax this restriction on the possible
extensions somewhat, and they are defined as follows: A minimal subshift $\Omega$ is
called an \textit{Arnoux-Rauzy subshift} if $p_\Omega(1) = m$ and for every $n \in \Z^+$,
there is a unique right-special word in $\mathcal{W}_\Omega(n)$ and a unique left-special
word in $\mathcal{W}_\Omega(n)$, both having exactly $m$ extensions. This implies in
particular that $p_\Omega(n) = (m-1)n + 1$. Arnoux-Rauzy subshifts over $\mathcal{A}_2$
are exactly the Sturmian subshifts.

A minimal subshift $\Omega$ is called \textit{episturmian} if $\mathcal{W}_\Omega$ is
closed under reversal (i.e., for every $w = w_1 \ldots w_n \in \mathcal{W}_\Omega$, we
have $w^R = w_n \ldots w_1 \in \mathcal{W}_\Omega$) and for every $n \in \Z^+$, there is
exactly one right-special word in $\mathcal{W}_\Omega(n)$.

\begin{prop}
Every Arnoux-Rauzy subshift is episturmian and every epi\-sturmian subshift is strictly
ergodic.
\end{prop}

See \cite{djp,juspir2,rz,wz} for these results and more information on Arnoux-Rauzy and
episturmian subshifts.

\subsubsection{Codings of Interval Exchange Transformations}

Interval exchange transformations are defined as follows. Given a probability vector
$\lambda = (\lambda_1,\ldots,\lambda_m)$ with $\lambda_i > 0$ for $1 \le i \le m$, let
$\mu_0 = 0$, $\mu_i = \sum_{j = 1}^i \lambda_j$, and $I_i = [\mu_{i-1},\mu_i)$. Let $\tau
\in S_m$, the symmetric group. Then $\lambda^\tau = (\lambda_{\tau^{-1}(1)}, \ldots,
\lambda_{\tau^{-1}(m)})$ is also a probability vector and we can form the corresponding
$\mu_i^\tau$ and $I_i^\tau$. Denote the unit interval $[0,1)$ by $I$. The
$(\lambda,\tau)$ interval exchange transformation is then defined by
$$
T : I \to I, \; \; T(x) = x - \mu_{i-1} + \mu_{\tau(i) - 1}^\tau \text{ for } x \in I_i,
\; 1 \le i \le m.
$$
It exchanges the intervals $I_i$ according to the permutation $\tau$.

The transformation $T$ is invertible and its inverse is given by the
$(\lambda^\tau,\tau^{-1})$ interval exchange transformation.

The symbolic coding of $x \in I$ is $\omega_n(x) = i$ if $T^n(x) \in I_i$. This induces a
subshift over the alphabet $\mathcal{A} = \{1,\ldots,m\}$: $\Omega_{\lambda,\tau} =
\overline{ \{ \omega(x) : x \in I \} }$. Every Sturmian subshift can be described by the
exchange of two intervals.

Keane \cite{Ke1} proved that if the orbits of the discontinuities $\mu_i$ of $T$ are all
infinite and pairwise distinct, then $T$ is minimal. In this case, the coding is
one-to-one and the subshift is minimal and aperiodic.  This holds in particular if $\tau$
is irreducible and $\lambda$ is irrational. Here, $\tau$ is called irreducible if
$\tau(\{1,\ldots,k\}) \not= (\{1,\ldots,k\})$ for every $k < m$ and $\lambda$ is called
irrational if the $\lambda_i$ are rationally independent.

Keane also conjectured that all minimal interval exchange transformations give rise to a
uniquely ergodic system. This was disproved by Keynes and Newton \cite{KN} using five
intervals, and then by Keane \cite{Ke2} using four intervals (the smallest possible
number). The conjecture was therefore modified in \cite{Ke2} and then ultimately proven
by Masur \cite{m}, Veech \cite{Vee1}, and Boshernitzan \cite{Bosh1}: For every
irreducible $\tau \in S_m$ and for Lebesgue almost every $\lambda$, the subshift
$\Omega_{\lambda,\tau}$ is uniquely ergodic.

\subsubsection{Substitution Sequences}
All the previous examples were generalizations of Sturmian sequences. We now discuss a
class of examples that generalize a certain aspect of the Fibonacci sequence
$$
s_n = \chi_{[1-\theta,0)}(n \theta) , \quad \theta = \frac{\sqrt{5}-1}{2}.
$$
We saw above (see \eqref{fibwkrec} and its discussion) that this sequence, restricted to
the right half line, is obtained by iterating the map \eqref{fibsubst}. That is,
$$
1 \mapsto 10 \mapsto 101 \mapsto 10110 \mapsto 10110101 \mapsto \cdots
$$
has $(s_n)_{n \ge 1}$ as its limit. In other words, $(s_n)_{n \ge 1}$ is invariant under
the substitution rule \eqref{fibsubst}.

\begin{defi}[substitution]
Denote the set of words over the alphabet $\mathcal{A}$ by $\mathcal{A}^*$. A map $S :
\mathcal{A} \to \mathcal{A}^*$ is called a substitution. The naturally induced maps on
$\mathcal{A}^*$ and $\mathcal{A}^{\Z^+}$ are denoted by $S$ as well.
\end{defi}

\noindent\textbf{Examples.} (a) Fibonacci: $1 \mapsto 10 , \; 0 \mapsto 1$ \\
(b) Thue-Morse: $1 \mapsto 10 , \; 0 \mapsto 01$\\
(c) Period doubling: $1 \mapsto 10 , \; 0 \mapsto 11$\\
(d) Rudin-Shapiro: $1 \mapsto 12$, $2 \mapsto 13$, $3 \mapsto 42$, $4 \mapsto 43$

\begin{defi}[substitution sequence]
Let $S$ be a substitution. A sequence $s \in \mathcal{A}^{\Z^+}$ is called a substitution
sequence if it is a fixed point of $S$.
\end{defi}

If $S(a)$ begins with the symbol $a$ and has length at least two, it follows that
$|S^n(a)| \to \infty$ as $n \to \infty$ and $S^n(a)$ has $S^{n-1}(a)$ as a prefix. Thus,
the limit of $S^n(a)$ as $n \to \infty$ defines a substitution sequence $s$. In the
examples above, we obtain the following substitution sequences.\\[1mm]
(a) Fibonacci: $s_\mathrm{F} = 1011010110110 \ldots$ \\
(b) Thue-Morse: $s^{(1)}_\mathrm{TM} = 100101100110 \ldots$ and $s^{(0)}_\mathrm{TM} = 0110100110010110 \ldots$\\
(c) Period doubling: $s_\mathrm{PD} = 101110101011101110 \ldots$\\
(d) Rudin-Shapiro: $s^{(1)}_\mathrm{RS} = 1213124212134313 \ldots$ and
$s^{(4)}_\mathrm{RS} = 4342431343421242 \ldots$

\medskip

We want to associate a subshift $\Omega_s$ with a substitution sequence $s$. Since the
iteration of $S$ on a suitable symbol $a$ naturally defines a one-sided sequence $s$, we
have to alter the definition of $\Omega_s$ used above slightly. One possible way is to
extend $s$ to a two-sided sequence $\tilde s$ arbitrarily and then define
$$
\Omega_s = \{ \omega \in \mathcal{A}^\Z : \omega = T^{n_j} \tilde s \text{ for some
sequence } n_j \to \infty \}.
$$
A different way is to define $\Omega_s$ to be the set of all $\omega$'s with
$\mathcal{W}_\omega \subseteq \mathcal{W}_s$. Below we will restrict our attention to
so-called primitive substitutions and for them, these two definitions are equivalent.

To ensure that $\Omega_s$ is strictly ergodic, we need to impose some conditions on $S$.
A very popular sufficient condition is primitivity.

\begin{defi}
A substitution $S$ is called primitive if there is $k \in \Z^+$ such that for every pair
$a,b \in \mathcal{A}$, $S^k(a)$ contains the symbol $b$.
\end{defi}

It is easy to check that our four main examples are primitive. Moreover, if $S$ is
primitive, then every power of $S$ is primitive. Thus, even if $S(a)$ does not begin with
$a$ for any symbol $a \in \mathcal{A}$, we may replace $S$ by a suitable $S^m$ and then
find such an $a$, which in turn yields a substitution sequence associated with $S^m$ by
iteration.

\begin{theorem}
Suppose $S$ is primitive and $s$ is an associated substitution sequence. Then, $s$ is
linearly recurrent. Consequently, $\Omega_s$ is strictly ergodic.
\end{theorem}

See \cite{dz2,duhosk}. Linear recurrence clearly also implies that $p_s(n) = O(n)$. Fixed
points of non-primitive substitutions may have quadratic complexity. However, there are
non-primitive substitutions that have fixed points which are linearly recurrent and hence
define strictly ergodic subshifts; see \cite{dl5} for a characterization of linearly
recurrent substitution generated subshifts.

\subsubsection{Subshifts with Positive Topological Entropy}

All the examples discussed so far have linearly bounded complexity. One may wonder if
strict ergodicity places an upper bound on the growth of the complexity function. Here we
want to mention the existence of strictly ergodic subshifts that have a very fast growing
complexity function. In fact, it is possible to have growth that is arbitrarily close to
the maximum possible one on a logarithmic scale.

Given a sequence $s$ over an alphabet $\mathcal{A}$, $|\mathcal{A}| \ge 2$, its
(\textit{topological}) \textit{entropy} is given by
$$
h_s = \lim_{n \to \infty} \frac1n \log p_s(n).
$$
The existence of the limit follows from the fact that $n \mapsto \log p_s(n)$ is
subadditive. Moreover,
$$
0 \le h_s \le \log |\mathcal{A}|,
$$
where $|\cdot|$ denotes cardinality.

The following was shown by Hahn and Katznelson \cite{hk}:

\begin{theorem}
{\rm (a)} If $s$ is a uniformly recurrent sequence over the alphabet $\mathcal{A}$, then
$h_s < \log |\mathcal{A}|$.\\
{\rm (b)} For every $\delta \in (0,1)$, there are alphabets $\mathcal{A}^{(j)}$ and
sequences $s^{(j)}$ over $\mathcal{A}^{(j)}$, $j \in \Z^+$, such that $|
\mathcal{A}^{(j)}| \to \infty$ as $j \to \infty$, $h_{s^{(j)}} \ge \log \left[
|\mathcal{A}^{(j)}| (1-\delta) \right]$, and every $\Omega_{s^{(j)}}$ is strictly
ergodic.
\end{theorem}

\section{Associated Schr\"odinger Operators and Basic Results}\label{Sec3}

In this section we associate Schr\"odinger operators with a subshift $\Omega$ and a
sampling function $f$ mapping $\Omega$ to the real numbers. In subsequent sections we
will study spectral and dynamical properties of these operators.

Let $\Omega$ be a strictly ergodic subshift with invariant measure $\mu$ and let $f :
\Omega \to \R$ be continuous. Then, for every $\omega \in \Omega$, we define a potential
$V_\omega : \Z \to \R$ by $V_\omega(n) = f(T^n \omega)$ and a bounded operator $H_\omega$
acting on $\ell^2(\Z)$ by
$$
[H_\omega \psi](n) = \psi(n+1) + \psi(n-1) + V_\omega(n) \psi(n).
$$

\noindent\textbf{Example.} The most common choice for $f$ is $f(\omega) = g(\omega_0)$
with some $g : \mathcal{A} \to \R$. This is a special case of a locally constant function
that is completely determined by the values of $\omega_n$ for $n$'s from a finite window
around the origin, that is, $f$ is called \textit{locally constant} if it is of the form
$f(\omega) = h(\omega_{-M} \ldots \omega_N)$ for suitable integers $M,N \ge 0$ and $h :
\mathcal{A}^{N+M+1} \to \R$. Clearly, every locally constant $f$ is continuous.

\medskip

The family $\{ H_\omega \}_{\omega \in \Omega}$ is an ergodic family of discrete
one-dimensional Schr\"odinger operators in the sense of Carmona and Lacroix \cite{cl}. By
the general theory it follows that the spectrum and the spectral type of $H_\omega$ are
$\mu$-almost surely $\omega$-independent \cite[Sect.~V.2]{cl}:

\begin{theorem}
There exist sets $\Omega_0 \subseteq \Omega$, $\Sigma, \Sigma_\mathrm{pp},
\Sigma_\mathrm{sc}, \Sigma_\mathrm{ac} \subseteq \R$ such that $\mu(\Omega_0) = 1$ and
\begin{align}
\label{specconst} \sigma(H_\omega) & = \Sigma \\
\sigma_\mathrm{pp} (H_\omega) & = \Sigma_\mathrm{pp} \\
\sigma_\mathrm{sc} (H_\omega) & = \Sigma_\mathrm{sc} \\
\label{acspecconst} \sigma_\mathrm{ac} (H_\omega) & = \Sigma_\mathrm{ac}
\end{align}
for every $\omega \in \Omega_0$.
\end{theorem}

Here, $\sigma(H),\sigma_\mathrm{pp}(H),\sigma_\mathrm{sc}(H),\sigma_\mathrm{ac}(H)$
denote the spectrum, the closure of the set of eigenvalues, the singular continuous
spectrum and the absolutely continuous spectrum of the operator $H$, respectively.

Since $\Omega$ is minimal and $f$ is continuous, a simple argument involving strong
approximation shows that \eqref{specconst} even holds everywhere, rather than almost
everywhere:

\begin{theorem}
For every $\omega \in \Omega$, $\sigma(H_\omega) = \Sigma$.
\end{theorem}

\begin{proof}
By symmetry it suffices to show that for every pair $\omega_1,\omega_2 \in \Omega$,
$\sigma(H_{\omega_1}) \subseteq \sigma(H_{\omega_2})$. Due to minimality, there exists a
sequence $(n_j)_{j \ge 1}$ such that $T^{n_j} \omega_2 \to \omega_1$ as $j \to \infty$.
By continuity of $f$, $H_{T^{n_j} \omega_2}$ converges strongly to $H_{\omega_1}$ as $j
\to \infty$. Thus,
$$
\sigma(H_{\omega_1}) \subseteq \overline{\bigcup_{j \ge 1} \sigma(H_{T^{n_j} \omega_2})}
= \sigma(H_{\omega_2}).
$$
Here, the first step follows by strong convergence and the second step is a consequence
of the fact that each of the operators $H_{T^{n_j} \omega_2}$ is unitarily equivalent to
$H_{\omega_2}$ and hence has the same spectrum.
\end{proof}

Far more subtle is the result that \eqref{acspecconst} also holds everywhere:

\begin{theorem}\label{lskthm}
For every $\omega \in \Omega$, $\sigma_\mathrm{ac}(H_\omega) = \Sigma_\mathrm{ac}$.
\end{theorem}

For strictly ergodic models, such as the ones considered here, there are two proofs of
Theorem~\ref{lskthm} in the literature. It was shown, based on unique ergodicity, by
Kotani in \cite{k4}. A proof based on minimality was given by Last and Simon in
\cite{ls}.

\medskip

A map $A \in C(\Omega,\mathrm{SL}(2,\R))$ induces an $\mathrm{SL}(2,\R)$-cocycle over $T$
as follows:
$$
\tilde A : \Omega \times \R^2 \to \Omega \times \R^2, \; (\omega,v) \mapsto (T \omega,
A(\omega) v ).
$$
Note that when we iterate this map $n$ times, we get
$$
\tilde A^n(\omega,v) = (T^n \omega , A_n(\omega) v),
$$
where $A_n(\omega) = A(T^{n-1}\omega) \cdots A(\omega)$. We are interested in the
asymptotic behavior of the norm of $A_n(\omega)$ as $n \to \infty$. The multiplicative
ergodic theorem ensures the existence of $\gamma_A \ge 0$, called the \textit{Lyapunov
exponent}, such that
\begin{align*}
\gamma_A & = \lim_{n \to \infty} \frac1n \int \log \| A_n(\omega)\| \,
d\mu(\omega) \\
& = \inf_{n \ge 1} \frac1n \int \log \| A_n(\omega) \| \, d\mu(\omega) \\
& = \lim_{n \to \infty} \frac1n \log \| A_n(\omega) \| \quad \text{ for $\mu$-a.e. }
\omega.
\end{align*}

In the study of the operators $H_\omega$ the following cocycles are relevant:
$$
A^{f,E}(\omega) = \left( \begin{array}{cr} E - f(\omega) & -1 \\ 1 & 0 \end{array}
\right),
$$
where $f$ is as above and $E$ is a real number, called the \textit{energy}. We regard $f$
as fixed and write $\gamma(E)$ instead of $\gamma_{A^{f,E}}$ to indicate that our main
interest is in the mapping $E \mapsto \gamma(E)$. Let
$$
\mathcal{Z} = \{ E \in \R : \gamma(E) = 0 \}.
$$
Note that we leave the dependence on $\Omega$ and $f$ implicit.

These cocycles are important in the study of $H_\omega$ because $A_n^{f,E}$ is the
transfer matrix for the associated difference equation. That is, a sequence $u$ solves
\begin{equation}\label{eve}
u(n+1) + u(n-1) + V_\omega(n) u(n) = E u(n)
\end{equation}
if and only if it solves
$$
\left( \begin{array}{c} u(n) \\ u(n-1) \end{array} \right) = A_n^{f,E} \left(
\begin{array}{c} u(0) \\ u(-1) \end{array} \right),
$$
as is readily verified.

\section{Absence of Absolutely Continuous Spectrum}\label{Sec5}

Let $\Omega$ be a strictly ergodic subshift and $f : \Omega \to \R$ locally constant. It
follows that the resulting potentials $V_\omega$ take on only finitely many values. In
this section we study the absolutely continuous spectrum of $H_\omega$, equal to
$\Sigma_\mathrm{ac}$ for every $\omega \in \Omega$ by Theorem~\ref{lskthm}. In 1982,
Kotani made one of the deepest and most celebrated contributions to the theory of ergodic
Schr\"odinger operators by showing that $\Sigma_\mathrm{ac}$ is completely determined by
the Lyapunov exponent, or rather the set $\mathcal{Z}$. Namely, his results, together
with earlier ones by Ishii and Patur, show that $\Sigma_\mathrm{ac}$ is given by the
essential closure of $\mathcal{Z}$. In 1989, Kotani found surprisingly general
consequences of his theory in the case of potentials taking on finitely many values. We
will review these results below.

\medskip

By assumption, the potentials $V_\omega$ take values in a fixed finite subset
$\mathcal{B}$ of $\R$. Thus, they can be regarded as elements of $\mathcal{B}^\Z$,
equipped with product topology. Let $\nu$ be the measure on $\mathcal{B}^\Z$ which is the
push-forward of $\mu$ under the mapping
$$
\Omega \to \mathcal{B}^\Z, \quad \omega \to V_\omega.
$$
Recall that the support of $\nu$, denoted by $\mathrm{supp} \, \nu$, is the complement of
the largest open set $U$ with $\nu(U) = 0$. Let
\begin{align*}
\left( \mathrm{supp} \, \nu \right)_+ & = \big\{ V |_{\Z^+_0} : V \in \mathrm{supp} \,
\nu \big\} \\
\left( \mathrm{supp} \, \nu \right)_- & = \big\{ V |_{\Z^-} : V \in \mathrm{supp} \, \nu
\big\},
\end{align*}
where $\Z^+_0 = \{ 0,1,2,\ldots \}$ and $\Z^- = \{ \ldots, -3, -2 , -1 \}$.

\begin{defi}
The measure $\nu$ is called deterministic if every $V_+ \in \left( \mathrm{supp} \, \nu
\right)_+$ comes from a unique $V \in \mathrm{supp} \, \nu$ and every $V_- \in \left(
\mathrm{supp} \, \nu \right)_-$ comes from a unique $V \in \mathrm{supp} \, \nu$.
\end{defi}

Consequently, if $\nu$ is deterministic, there is a bijection $C : \left( \mathrm{supp}
\, \nu \right)_- \to \left( \mathrm{supp} \, \nu \right)_+$ such that for every $V \in
\mathrm{supp} \, \nu$, $V |_{\Z^+_0} = C (V |_{\Z^-})$ and $V |_{\Z^-} = C^{-1} (V
|_{\Z^+_0})$.

\begin{defi}
The measure $\nu$ is called topologically deterministic if it is deterministic and the
map $C$ is a homeomorphism.
\end{defi}

Thus, when $\nu$ is topologically deterministic, we can continuously recover one half
line from the other for elements of $ \mathrm{supp} \, \nu$.

Let us only state the part of Kotani theory that is of immediate interest to us here:

\begin{theorem}\label{kotthm1}
{\rm (a)} If $\mathcal{Z}$ has zero Lebesgue measure, then $\Sigma_\mathrm{ac}$ is
empty. \\
{\rm (b)} If $\mathcal{Z}$ has positive Lebesgue measure, then $\nu$ is topologically
deterministic.
\end{theorem}

This theorem holds in greater generality; see \cite{k2,k3,k4,s1}. The underlying
dynamical system $(\Omega,T,\mu)$ is only required to be measurable and ergodic and the
set $\mathcal{B}$ can be any compact subset of $\R$. Part~(a) is a particular consequence
of the Ishii-Kotani-Pastur identity
$$
\Sigma_\mathrm{ac} = \overline{\mathcal{Z}}^\mathrm{ess},
$$
where the essential closure of a set $S \subseteq \R$ is given by
$$
\overline{S}^\mathrm{ess} = \{ E \in \R : \mathrm{Leb} \left( (E - \varepsilon,E +
\varepsilon) \cap S \right) > 0 \text{ for every } \varepsilon > 0 \}.
$$

The following result was proven by Kotani in 1989 \cite{k3}. Here it is crucial that the
set $\mathcal{B}$ is finite.

\begin{theorem}\label{kotthm2}
If $\nu$ is topologically deterministic, then $\mathrm{supp} \, \nu$ is finite.
Consequently, all potentials in $\mathrm{supp} \, \nu$ are periodic.
\end{theorem}

Combining Theorems~\ref{kotthm1} and \ref{kotthm2} we arrive at the following corollary.

\begin{coro}\label{kotcoro}
If $\Omega$ and $f$ are such that $\mathrm{supp} \, \nu$ contains an aperiodic element,
then $\mathcal{Z}$ has zero Lebesgue measure and $\Sigma_\mathrm{ac}$ is empty.
\end{coro}

Note that by minimality, the existence of one aperiodic element is equivalent to all
elements being aperiodic. This completely settles the issue of existence/purity of
absolutely continuous spectrum. In the periodic case, the spectrum of $H_\omega$ is
purely absolutely continuous for every $\omega \in \Omega$, and in the aperiodic case,
the spectrum of $H_\omega$ is purely singular for every $\omega \in \Omega$.

\section{Zero-Measure Spectrum}\label{Sec6}

Suppose throughout this section that $\Omega$ is strictly ergodic, $f : \Omega \to \R$ is
locally constant, and the resulting potentials $V_\omega$ are aperiodic.\footnote{Even
when we make explicit assumptions on $\Omega$ and $f$, aperiodicity of the potentials
will always be assumed implicitly; for example, in Theorem~\ref{primsubstzms} and
Corollary~\ref{lenzcoro}.} This section deals with the Lebesgue measure of the set
$\Sigma$, which is the common spectrum of the operators $H_\omega$, $\omega \in \Omega$.
It is widely expected that $\Sigma$ always has zero Lebesgue measure. This is supported
by positive results for large classes of subshifts and functions. We present two
approaches to zero-measure spectrum, one based on trace map dynamics and sub-exponential
upper bounds for $\|A^{f,E}_n(\omega)\|$ for energies in the spectrum, and another one
based on uniform convergence of $\frac1n \log \| A^{f,E}_n(\omega) \|$ to $\gamma(E)$ for
all energies. Both approaches have in common that they establish the identity
\begin{equation}\label{sigmaz}
\Sigma = \mathcal{Z}.
\end{equation}
Zero-measure spectrum then follows immediately from Corollary~\ref{kotcoro}.

\subsection{Trace Map Dynamics}

Zero-measure spectrum follows once one proves
\begin{equation}\label{sinz}
\Sigma \subseteq \mathcal{Z}.
\end{equation}
Note, however, that the Lyapunov exponent is always positive away from the spectrum.
Thus, $\mathcal{Z} \subseteq \Sigma$, and \eqref{sinz} is in fact equivalent to
\eqref{sigmaz}.

A trace map is a dynamical system that may be associated with a family $\{H_\omega\}$
under suitable circumstances. It is given by the iteration of a map $T : \R^k \to \R^k$.
Iteration of this map on some energy-dependent initial vector, $v_E$, will then describe
the evolution of a certain sequence of transfer matrix traces. Typically, these iterates
will diverge rather quickly. The stable set, $B_\infty$, is defined to be the set of
energies for which $T^n v_E$ does not diverge quickly. The inclusion \eqref{sinz} is then
established in a two-step procedure:
\begin{equation}\label{sinbinz}
\Sigma \subseteq B_\infty \subseteq \mathcal{Z}.
\end{equation}
Again, by the remark above, this establishes equality and hence
$$
\Sigma = B_\infty = \mathcal{Z}.
$$

For the sake of clarity of the main ideas, we first discuss the trace-map approach for
the Fibonacci subshift $\Omega_\mathrm{F}$ and $f : \Omega_\mathrm{F} \to \R$ given by
$f(\omega) = g(\omega_0)$, $g(0) = 0$, $g(1) = \lambda > 0$. See \cite{c,kkt,oprss,s5,s6}
for the original literature concerning this special case.

Given the partition result, Theorem~\ref{partition}, and the recursion \eqref{fibwkrec},
it is natural to decompose transfer matrix products into factors of the form $M_k(E)$,
where
\begin{equation}\label{mkrecbasis}
M_{-1}(E) = \left( \begin{array}{cr} 1 & -\lambda \\ 0 & 1 \end{array} \right), \quad
M_0(E) = \left( \begin{array}{cr} E & -1 \\ 1 & 0 \end{array} \right)
\end{equation}
and
\begin{equation}\label{fibmkrec}
M_{k+1}(E) = M_{k-1}(E) M_k(E), \; \text{ for } k \ge 0.
\end{equation}

\begin{prop}
Let $x_k = x_k(E) = \frac12 \mathrm{Tr} M_k(E)$. Then,
\begin{equation}\label{fibtracemap}
x_{k+2} = 2 x_{k+1} x_k - x_{k-1} \quad \text{ for } k \in \Z^+_0
\end{equation}
and
\begin{equation}\label{fibinv}
x_{k+1}^2 + x_k^2 + x_{k-1}^2 - 2 x_{k+1} x_k x_{k-1} = 1 + \frac{\lambda^2}{4} \quad
\text{ for } k \in \Z^+_0.
\end{equation}
\end{prop}

\begin{proof}
The recursion \eqref{fibtracemap} follows readily from \eqref{fibmkrec}. Using
\eqref{fibtracemap}, one checks that the left-hand side of \eqref{fibinv} is independent
of $k$. Evaluation for $k = 0$ then yields the right-hand side. See \cite{kkt,oprss,s5}
for more details.
\end{proof}

The recursion \eqref{fibtracemap} is called the \textit{Fibonacci trace map}. The $x_k$'s
may be obtained by the iteration of the map $T : \R^3 \to \R^3$, $(x, y, z) \mapsto (xy -
z, x, y)$ on the initial vector $((E-\lambda)/2,E/2,1)$.

\begin{prop}\label{fibdicho}
The sequence $(x_k)_{k \ge -1}$ is unbounded if and only if
\begin{equation}\label{escapecond}
|x_{k_0-1}| \le 1, \quad |x_{k_0}| > 1, \quad |x_{k_0+1}| > 1
\end{equation}
for some $k_0 \ge 0$. In this case, the $k_0$ is unique, and we have
\begin{equation}\label{expfib}
|x_{k+2}| > |x_{k+1} x_k| > 1 \quad \text{ for } k \ge k_0
\end{equation}
and
\begin{equation}\label{expfib2}
|x_k| > C^{F_{k-k_0}} \quad \text{ for } k \ge k_0
\end{equation}
and some $C > 1$. If $(x_k)_{k \ge -1}$ is bounded, then
\begin{equation}\label{univfibbound}
|x_k| \le 1 + \frac{\lambda}{2} \quad \text{ for every } k.
\end{equation}
\end{prop}

\begin{proof}
Suppose first that \eqref{escapecond} holds for some $k_0 \ge 0$. Then, by
\eqref{fibtracemap},
$$
|x_{k_0+2}| \ge |x_{k_0+1}x_{k_0}| + ( |x_{k_0+1}x_{k_0}| - |x_{k_0-1}| ) >
|x_{k_0+1}x_{k_0}| > 1.
$$
By induction, we get \eqref{expfib}, and also that the $k_0$ is unique. Taking $\log$'s,
we see that $\log |x_k|$ grows faster than a Fibonacci sequence for $k \ge k_0$, which
gives \eqref{expfib2}.

Conversely, suppose that \eqref{escapecond} fails for every $k_0 \ge 0$. Consider a value
of $k$ for which $|x_k| > 1$. Since $x_{-1} = 1$, it follows that $|x_{k-1}| \le 1$ and
$|x_{k+1}| \le 1$. Thus, the invariant \eqref{fibinv} shows that
\begin{align*}
|x_k| & \le |x_{k+1}x_{k-1}| + \left( |x_{k+1}x_{k-1}|^2 - x_{k+1}^2 - x_{k-1}^2 + 1 +
\frac{\lambda^2}{4} \right)^{1/2} \\
 & = |x_{k+1}x_{k-1}| + \left( (1 - x_{k+1}^2)(1 - x_{k-1}^2) +
\frac{\lambda^2}{4} \right)^{1/2},
\end{align*}
which implies that the sequence $(x_k)_{k \ge -1}$ is bounded and obeys
\eqref{univfibbound}.
\end{proof}

The dichotomy described in Proposition~\ref{fibdicho} motivates the following definition:
\begin{equation}\label{binfdef}
B_\infty = \left\{ E \in \R : |x_k| \le 1 + \frac{\lambda}{2} \text{ for every } k
\right\}.
\end{equation}
This set provides the link between the spectrum and the set of energies for which the
Lyapunov exponent vanishes.

\begin{theorem}\label{fibzms}
Let $\Omega = \Omega_\mathrm{F}$ be the Fibonacci subshift and let $f : \Omega \to \R$ be
given by $f(\omega) = g(\omega_0)$, $g(0) = 0$, $g(1) = \lambda > 0$. Then, $\Sigma =
B_\infty = \mathcal{Z}$ and $\Sigma$ has zero Lebesgue measure.
\end{theorem}

\begin{proof}
We show the two inclusions in \eqref{sinbinz}. Let $\sigma_k = \{ E : |x_k| \le 1 \}$. On
the one hand, $\sigma_k$ is the spectrum of an $F_k$-periodic Schr\"odinger operator
$H_k$. It is not hard to see that $H_k \to H$ strongly, where $H$ is the Schr\"odinger
operator with potential $V(n) = \lambda \chi_{[1-\theta,1)}(n\theta)$ and hence $\Sigma$
is contained in the closure of $\bigcup_{\tilde k \ge k} \sigma_{\tilde k}$ for every
$k$. On the other hand, Proposition~\ref{fibdicho} shows that $\sigma_{k+1} \cup
\sigma_{k+2} \subseteq \sigma_k \cup \sigma_{k+1}$ and $B_\infty = \bigcap_{k} \sigma_k
\cup \sigma_{k+1}$. Thus,
$$
\Sigma \subseteq \bigcap_{k} \overline{ \bigcup_{\tilde k \ge k} \sigma_{\tilde k} } =
\bigcap_{k} \sigma_k \cup \sigma_{k+1} = B_\infty.
$$
This is the first inclusion in \eqref{sinbinz}. The second inclusion follows once we can
show that for every $E \in B_\infty$, we have that $\log \| M_n \| \lesssim n$, where the
implicit constant depends only on $\lambda$. From the matrix recursion \eqref{fibmkrec}
and the Cayley-Hamilton Theorem, we obtain
$$
M_{k+1} = M_{k-1} M_k^2 M_k^{-1} = M_{k-1} (2 x_k M_k - \mathrm{Id}) M_k^{-1} = 2 x_k
M_{k-1} - M_{k-2}^{-1}.
$$
If $E \in B_\infty$, then $2|x_k| \le 2 + \lambda$, and hence we obtain by induction that
$\|M_k\| \le C^k$. Combined with the partition result, Theorem~\ref{partition}, this
yields the claim since the $F_k$ grow exponentially.
\end{proof}

The same strategy works in the Sturmian case, as shown by Bellissard et al.\ \cite{bist},
although the analysis is technically more involved. Because of \eqref{wkrec}, we now
consider instead of \eqref{fibmkrec} the matrices defined by the recursion
$$
M_{k+1}(E) = M_{k-1}(E) M_k(E)^{a_{k+1}},
$$
where the $a_k$'s are the coefficients in the continued fraction expansion \eqref{cfexp}
of $\theta$. This recursion again gives rise to a trace map for $x_k = \frac12
\mathrm{Tr} M_k(E)$ which involves Chebyshev polynomials. These traces obey the invariant
\eqref{fibinv} and the exact analogue of Proposition~\ref{fibdicho} holds. After these
properties are established, the proof may be completed as above. Namely, $B_\infty$ is
again defined by \eqref{binfdef} and the same line of reasoning yields the two inclusions
in \eqref{sinbinz}. We can therefore state the following result:

\begin{theorem}\label{sturmzms}
Let $\Omega$ be a Sturmian subshift with irrational slope $\theta \in (0,1)$ and let $f :
\Omega \to \R$ be given by $f(\omega) = g(\omega_0)$, $g(0) = 0$, $g(1) = \lambda > 0$.
Then, $\Sigma = B_\infty = \mathcal{Z}$ and $\Sigma$ has zero Lebesgue measure.
\end{theorem}

We see that every operator family associated with a Sturmian subshift admits a trace map
and an analysis of this dynamical system allows one to prove the zero-measure property.

Another class of operators for which a trace map always exists and may be used to prove
zero-measure spectrum is given by those that are generated by a primitive substitution.
The existence of a trace map is even more natural in this case and not hard to verify;
see, for example, \cite{ap,ab,abg,kn,pww} for general results and \cite{bgj,br,rb} for
trace maps with an invariant. However, its analysis is more involved and has been
completed only in 2002 by Liu et al.\ \cite{ltww}, after a number of earlier works had
established partial results \cite{b2,bbg1,bg1}. Bellissard et al., on the other hand, had
proved their Sturmian result already in 1989 -- shortly after Kotani made his crucial
observation leading to Corollary~\ref{kotcoro}.

\begin{theorem}\label{primsubstzms}
Let $\Omega$ be a subshift generated by a primitive substitution $S : \mathcal{A} \to
\mathcal{A}^*$ and let $f : \Omega \to \R$ be given by $f(\omega) = g(\omega_0)$ for some
function $g : \mathcal{A} \to \R$. Then, the associated trace map admits a stable set,
$B_\infty$, for which we have $\Sigma = B_\infty = \mathcal{Z}$. Consequently, $\Sigma$
has zero Lebesgue measure.
\end{theorem}

\subsection{Uniform Hyperbolicity}

Recall that the Lyapunov exponent associated with the Schr\"odinger cocycle $A^{f,E}$
obeys
\begin{equation}\label{leconv}
\gamma(E) = \lim_{n \to \infty} \frac1n \log \| A^{f,E}(T^{n-1} \omega) \cdots
A^{f,E}(\omega) \|
\end{equation}
for $\mu$-almost every $\omega \in \Omega$.

\begin{defi}[uniformity]
The cocycle $A^{f,E}$ is called uniform if the convergence in \eqref{leconv} holds for
every $\omega \in \Omega$ and is uniform in $\omega$. It is called uniformly hyperbolic
if it is uniform and $\gamma(E) > 0$. Define
$$
\mathcal{U} = \{ E \in \R : A^{f,E} \text{ is uniformly hyperbolic } \}.
$$
\end{defi}

Uniform hyperbolicity of $A^{f,E}$ is equivalent to $E$ belonging to the resolvent set as
shown by Lenz \cite{len2} (see also Johnson \cite{john}):

\begin{theorem}\label{lenzthm}
$\R \setminus \Sigma = \mathcal{U}$.
\end{theorem}

Recall that we assumed at the beginning of this section that the potentials $V_\omega$
are aperiodic. Thus, combining Corollary~\ref{kotcoro} and Theorem~\ref{lenzthm}, we
arrive at the following corollary.

\begin{coro}\label{lenzcoro}
If $A^{f,E}$ is uniform for every $E \in \R$, then $\Sigma = \mathcal{Z}$ and $\Sigma$
has zero Lebesgue measure.
\end{coro}

We thus seek a sufficient condition on $\Omega$ and $f$ such that $A^{f,E}$ is uniform
for every $E \in \R$, which holds for as many cases of interest as possible. Such a
condition was recently found by Damanik and Lenz in \cite{dl7}. In \cite{dl8} it was then
shown by the same authors that this condition holds for the majority of the models
discussed in Section~\ref{Sec2}.

\begin{defi}[condition {\rm (B)}]
Let $\Omega$ be a strictly ergodic subshift with unique $T$-invariant measure $\mu$. It
satisfies the Boshernitzan condition {\rm (B)} if
\begin{equation}\label{boshcond}
\limsup_{n \to \infty} \left( \min_{w \in \mathcal{W}_\Omega(n)} n \cdot \mu \left( [w]
\right) \right)
> 0.
\end{equation}
\end{defi}

\noindent\textbf{Remarks.} (a) $[w]$ denotes the cylinder set
$$
[w] = \{ \omega \in \Omega : \omega_1 \ldots \omega_{|w|} = w \}.
$$
(b) It suffices to assume that $\Omega$ is minimal and there exists some $T$-invariant
measure $\mu$ with \eqref{boshcond}. Then, $\Omega$ is necessarily uniquely ergodic.\\
(c) The condition \eqref{boshcond} was introduced by Boshernitzan in \cite{Bosh2}. His
main purpose was to exhibit a useful sufficient condition for unique ergodicity. The
criterion proved to be particularly useful in the context of interval exchange
transformations \cite{Bosh1,Vee2}, where unique ergodicity holds almost always, but not
always \cite{Ke2,m,Vee1}.

\medskip

It was shown by Damanik and Lenz that condition (B) implies uniformity for all energies
and hence zero-measure spectrum \cite{dl7}.

\begin{theorem}\label{boshthmdl}
If $\Omega$ satisfies {\rm (B)}, then $A^{f,E}$ is uniform for every $E \in \R$.
\end{theorem}

\noindent\textbf{Remarks.} (a) The Boshernitzan condition holds for almost all of the
subshifts discussed in Section~\ref{Sec2}. For example, it holds for every Sturmian
subshift, almost every subshift generated by a coding of a rotation with respect to a
two-interval partition, a dense set of subshifts associated with general codings of
rotations, almost every subshift associated with an interval exchange transformation,
almost every episturmian subshift, and every linearly recurrent subshift; see
\cite{dl8}.\footnote{Here, notions like ``dense'' or ``almost all'' are with respect to
the natural parameters of the class of models in question. We refer the reader to
\cite{dl8} for detailed statements of these applications of Theorem~\ref{boshthmdl}.} \\
(b) There was earlier work by Lenz who proved uniformity for all energies assuming a
stronger condition, called (PW) for positive weights \cite{len1}. Essentially, (PW)
requires \eqref{boshcond} with $\limsup$ replaced by $\liminf$. The condition (PW) holds
for all linearly recurrent subshifts but it fails, for example, for almost every Sturmian
subshift.\\
(c) Lenz in turn was preceded and inspired by Hof \cite{h} who proved uniform existence
of the Lyapunov exponents for Schr\"odinger operators associated with primitive
substitutions. Extensions of \cite{h} to linearly recurrent systems, including
higher-dimensional ones, were found by Damanik and Lenz \cite{dl11}.

\medskip

In particular, all zero-measure spectrum results obtained by the trace map approach also
follow from Theorem~\ref{boshthmdl}. Moreover, the applications of
Theorem~\ref{boshthmdl} cover operator families that are unlikely to be amenable to the
trace map approach. However, as we will see later, the trace map approach yields
additional information that is crucial in a study of detailed spectral and dynamical
properties. Thus, it is worthwhile to carry out an analysis of the trace map whenever
possible.

\subsection{The Hausdorff Dimension of the Spectrum}

Once one knows that the spectrum has zero Lebesgue measure it is a natural question if
anything can be said about its Hausdorff dimension. There is very important unpublished
work by Raymond \cite{r} who proved in the Fibonacci setting of Theorem~\ref{fibzms} that
the Hausdorff dimension is strictly smaller than one for $\lambda$ large enough ($\lambda
\ge 5$ is sufficient) and it converges to zero as $\lambda \to \infty$. Strictly positive
lower bounds for the Hausdorff dimension of the spectrum at all couplings $\lambda$
follow from the Hausdorff continuity results of \cite{d2,jl2}, to be discussed in
Section~\ref{Sec8}. Several aspects of Raymond's work were used and extended in a number
of papers \cite{d8,dt1,dt3,kkl,lw}. In particular, Liu and Wen carried out a detailed
analysis of the Hausdorff dimension of the spectrum in the general Sturmian case in the
spirit of Raymond's approach; see \cite{lw}.

\section{Absence of Point Spectrum}\label{Sec7}

We have seen that two of the three properties that are expected to hold in great
generality for the operators discussed in this paper hold either always (absence of
absolutely continuous spectrum) or almost always (zero-measure spectrum). In this section
we turn to the third property that is expected to be the rule---the absence of point
spectrum. As with zero-measure spectrum, no counterexamples are known and there are many
positive results that have been obtained by essentially two different methods. Both
methods rely on local symmetries of the potential. The existence of square-summable
eigenfunctions is excluded by showing that these local symmetries are reflected in the
solutions of the difference equation \eqref{eve}. Since there are exactly two types of
symmetries in one dimension, the effective criteria for absence of eigenvalues that
implement this general idea therefore rely on translation and reflection symmetries,
respectively. In the following, we explain these two methods and their range of
applicability.

\subsection{Local Repetitions}

In 1976, Gordon showed how to use the Cayley-Hamilton theorem to derive quantitative
solution estimates from local repetitions in the potential \cite{g2}. The first major
application of this observation was in the context of the almost Mathieu operator: For
super-critical coupling and Liouville frequencies, there is purely singular continuous
spectrum for all phases, as shown by Avron and Simon in 1982 \cite{avsim}.\footnote{As a
consequence, positive Lyapunov exponents do not in general imply spectral localization.}
The first application of direct relevance to this survey was obtained by Delyon and
Petritis \cite{dp1} in 1986 who proved absence of eigenvalues for certain codings of
rotations, including most Sturmian models. Further applications will be mentioned below.

Gordon's Lemma is a deterministic criterion and may be applied to a fixed potential $V :
\Z \to \R$. Analogous to the discussion in Section~\ref{Sec3}, we define transfer
matrices $A^E_n = T_{n-1} \cdots T_0$, where
$$
T_j = \left( \begin{array}{cr} E - V(j) & -1 \\ 1 & 0 \end{array} \right).
$$
Then, a sequence $u$ solves
\begin{equation}\label{evedet}
u(n+1) + u(n-1) + V(n) u(n) = E u(n)
\end{equation}
if and only if it solves $U(n) = A^E_n U(0)$, where
$$
U(j) = \left( \begin{array}{c} u(j) \\ u(j-1) \end{array} \right).
$$

\begin{lemma}\label{gordon2}
Suppose the potential $V$ obeys $V(m+p) = V(m)$, $0 \le m \le p-1$. Then,
$$
\max \left\{ \left\| U(2p) \right\|, \left\| U(p)  \right\| \right\} \ge \frac{1}{2 \max
\{ | \mathrm{Tr} A^E_p | , 1\}}  \left\| U(0) \right\|.
$$
\end{lemma}

\begin{proof}
This is immediate from the Cayley-Hamilton Theorem, applied to the matrix $A^E_p$ and the
vector $(u(0),u(-1))^T$.
\end{proof}

For obvious reasons, we call this criterion the \textit{two-block} Gordon Lemma. A slight
variation of the argument gives the following (\textit{three-block}) version of Gordon's
Lemma.

\begin{lemma}\label{gordon3}
Suppose the potential $V$ obeys $V(m+p) = V(m)$, $-p \le m \le p-1$. Then,
$$
\max \left\{ \left\| U(2p) \right\|, \left\| U(p)  \right\| ,
 \left\| U(-p) \right\| \right\} \ge \frac{1}{2}.
$$
\end{lemma}

\noindent\textbf{Remark.} The original criterion from \cite{g2} uses four blocks. For the
application to the almost Mathieu operator, this is sufficient; but for Sturmian models,
for example, the improvements above are indeed needed, as we will see below. The
two-block version can be found in S\"ut\H{o}'s paper \cite{s5} and the three-block
version was proved in \cite{dp1} by Delyon and Petritis.

\medskip

The two-block version is especially useful in situations where a trace map exists and we
have bounds on trace map orbits for energies in the spectrum. Note that the two-block
version gives a stronger conclusion. This will be crucial in the next section when we
discuss continuity properties of spectral measures with respect to Hausdorff measures in
the context of quantum dynamics.

\subsection{Palindromes}

Gordon-type criteria give quantitative estimates for solutions of \eqref{evedet} in the
sense that repetitions in the potential are reflected in solutions, albeit in a weak
sense. One would hope that local reflection symmetries in the potential give similar
information. Unfortunately, such a result has not been found yet. However, it is possible
to exclude square-summable solutions in this way by an indirect argument. Put slightly
simplified, \textit{if} a solution is square-summable, \textit{then} local reflection
symmetries are mirrored by solutions and these solution symmetries in turn prevent the
solution from being square-summable.

The original criterion for absence of eigenvalues in this context is due to Jitomirskaya
and Simon \cite{js} and it was developed in the context of the almost Mathieu operator to
prove, just as the result by Avron and Simon did, an unexpected occurrence of singular
continuous spectrum. An adaptation of the Jitomirskaya-Simon method to the subshift
context can be found in a paper by Hof et al.\ \cite{hks}. Let us state their
result:\footnote{There is also a half-line version, which is the palindrome analogue of
Lemma~\ref{gordon2}; see \cite{dgr}.}

\begin{lemma}\label{hkslem}
Let $V : \Z \to \R$ be given. There is a constant $B$, depending only on $\|V\|_\infty$,
with the following property: if there are $n_j \to \infty$ and $l_j$ with $B^{n_j}/l_j
\to 0$ as $j \to \infty$ such that $V$ is symmetric about $n_j$ on an interval of length
$l_j$ centered at $n_j$ for every $j$, then the Schr\"odinger operator $H$ with potential
$V$ has empty point spectrum.
\end{lemma}

\begin{proof}[Sketch of proof.]
Suppose that $V$ satisfies the assumptions of the lemma. Assume that $u$ is a
square-summable solution of \eqref{evedet}, normalized so that $\|u\|_2 = 1$. Fix some
$j$ and reflect $u$ about $n_j$. Call the reflected sequence $u^{(j)}$. Since the
potential is reflection-symmetric on an interval of length $l_j$ about $n_j$, the
Wronskian of $u$ and $u^{(j)}$ is constant on this interval. By $\|u\|_2 = 1$, it is
pointwise bounded in this interval by $2/l_j$. From this, it follows that $u$ and
$u^{(j)}$ are close (up to a sign) near $n_j$. Now apply transfer matrices and compare
$u$ and $u^{(j)}$ near zero. The assumption $B^{n_j}/l_j \to 0$ then implies that, for
$j$ large, $u$ and $u^{(j)}$ are very close near zero. In other words, $u$ is bounded
away from zero near $2n_j$ for all large $j$. This contradicts $u \in \ell^2(\Z)$.
\end{proof}

Thus, eigenvalues can be excluded if the potential contains infinitely many suitably
located palindromes. Here, a \textit{palindrome} is a word that is the same when read
backwards. Sequences obeying the assumption of Lemma~\ref{hkslem} are called
\textit{strongly palindromic} in \cite{hks}.

Hof et al.\ also prove the following general result for subshifts:

\begin{prop}\label{hksprop}
Suppose $\Omega$ is an aperiodic minimal subshift. If $\mathcal{W}_\Omega$ contains
infinitely many palindromes, then the set of strongly palindromic $\omega$'s in $\Omega$
is uncountably infinite.
\end{prop}

In any event, since the set $\mathcal{C}_\Omega = \{ \omega \in \Omega :
\sigma_\mathrm{pp}(H_\omega) = \emptyset \}$ is a $G_\delta$ set as shown by Simon
\cite{s2} (see also Choksi and Nadkarni \cite{chna} and Lenz and Stollmann \cite{lest}),
it is a dense $G_\delta$ set as soon as it is non-empty by minimality of $\Omega$ and
unitary equivalence of $H_\omega$ and $H_{T \omega}$.

Thus, when excluding eigenvalues we are interested in three kinds of results. We say that
eigenvalues are \textit{generically absent} if $\mathcal{C}_\Omega$ is a dense $G_\delta$
set. To prove generic absence of eigenvalues it suffices to treat one $\omega \in
\Omega$, as explained in the previous paragraph. Absence of eigenvalues holds
\textit{almost surely} if $\mu(\mathcal{C}_\Omega) =1$. To prove almost sure absence of
eigenvalues one only has to show $\mu(\mathcal{C}_\Omega) > 0$ by ergodicity and
$T$-invariance of $\mathcal{C}_\Omega$. Finally, absence of eigenvalues is said to hold
\textit{uniformly} if $\mathcal{C}_\Omega = \Omega$.

\subsection{Applications}

Let us now turn to applications of the two methods just described. We emphasize that
absence of eigenvalues is expected to hold in great generality and no counterexamples are
known.

As in Section~\ref{Sec6}, things are completely understood in the Sturmian case and
absence of eigenvalues holds uniformly.

\begin{theorem}\label{sturmapp}
Let $\Omega$ be a Sturmian subshift with irrational slope $\theta \in (0,1)$ and let $f :
\Omega \to \R$ be given by $f(\omega) = g(\omega_0)$, $g(0) = 0$, $g(1) = \lambda > 0$.
Then, $H_\omega$ has empty point spectrum for every $\omega \in \Omega$.
\end{theorem}

\begin{proof}[Sketch of proof.]
Given any $\lambda > 0$ and $\omega \in \Omega$, absence of point spectrum follows if
Lemma~\ref{gordon2} can be applied to $V_\omega$ for infinitely many values of $p$.
Considering only $p$'s of the form $q_k$, where the $q_k$'s are associated with $\theta$
via \eqref{qkdef}, the trace bounds established in Theorem~\ref{sturmzms} show that we
can focus our attention on the existence of infinitely many two-block structures aligned
at the origin. Using Theorem~\ref{partition}, a case-by-case analysis through the various
levels of the hierarchy detects these structures and completes the proof.
\end{proof}

\noindent\textbf{Remarks.} (a) For details, see \cite{dkl,dl1}. In fact, the argument
above has to be extended slightly for $\theta$'s with $\limsup a_k = 2$. To deal with
these exceptional cases, one also has to consider $p$'s of the form $q_k +
q_{k-1}$.\\
(b) Here is a list of earlier partial results for Sturmian models: Delyon and Petritis
proved absence of eigenvalues almost surely for every $\lambda > 0$ and Lebesgue almost
every $\theta$ \cite{dp1}. Their proof employs Lemma~\ref{gordon3}. Using
Lemma~\ref{gordon2}, S\"ut\H{o} proved absence of eigenvalues for $\lambda > 0$, $\theta
= (\sqrt{5}-1)/2$, and $\phi = 0$ \cite{s5}, and hence generic absence of eigenvalues in
the Fibonacci case. His proof and result were extended to all irrational $\theta$'s by
Bellissard et al.\ \cite{bist}.\footnote{They do not state the result explicitly in
\cite{bist}, but given their analysis of the trace map and the structure of the
potential, it follows as in \cite{s5}.} Hof et al.\ proved generic absence of eigenvalues
for every $\lambda > 0$ and every $\theta$ using Lemma~\ref{hkslem} \cite{hks}. Kaminaga
then showed an almost sure result for every $\lambda > 0$ and every $\theta$ \cite{k1}.
His proof is based on Lemma~\ref{gordon3} and refines the arguments of Delyon and Petritis.\\
(c) If most of the continued fraction coefficients are small, eigenvalues cannot be
excluded using a four-block Gordon Lemma. This applies in particular in the Fibonacci
case where $a_k \equiv 1$. The reason for this is that there simply are no four-block
structures in the potential. See \cite{dl10,dl9,juspir,van} for papers dealing with local
repetitions in Sturmian sequences.\\
(d) The palindrome method is very useful to prove generic results. However, it cannot be
used to prove almost sure or uniform results for linearly recurrent subshifts (e.g.,
subshifts generated by primitive substitutions). Namely, for these subshifts, the
strongly palindromic elements form a set of zero $\mu$-measure as shown by Damanik and
Zare \cite{dz2}.

\medskip

Let us now turn to subshifts generated by codings of rotations. The key papers were
mentioned above \cite{dp1,hks,k1}.

\begin{theorem}\label{corapp}
Suppose $\Omega$ is the subshift generated by a sequence of the form \eqref{codofrot}
with irrational $\theta \in (0,1)$, some $\phi \in [0,1)$, and a partition on the circle
into $l$ half-open intervals. Let $f : \Omega \to \R$ be given by $f(\omega) =
g(\omega_0)$ with some non-constant function $g$. Suppose that the continued fraction
coefficients of $\theta$ satisfy
\begin{equation}\label{dpcond}
\limsup_{k \to \infty} a_k \ge 2 l.
\end{equation}
Then, $H_\omega$ has empty point spectrum for $\mu$-almost every $\omega \in \Omega$.
\end{theorem}

\noindent\textbf{Remarks.} (a) For every $l \in \Z^+$, the condition \eqref{dpcond} holds
for Lebesgue almost every $\theta$. In fact, almost every $\theta$ has unbounded
continued
fraction coefficients; see \cite{khin}.\\
(b) The proof of Theorem~\ref{corapp}, given in \cite{dp1,k1}, is based on
Lemma~\ref{gordon3}.\\
(c) Hof et al.\ prove a generic result using Lemma~\ref{hkslem} for every $\theta$
provided that the partition of the circle has a certain symmetry property, which is
always satisfied in the case $l = 2$ \cite{hks}.\\
(d) It is possible to prove a result similar to Theorem~\ref{corapp} for a locally
constant $f$. In this case, the number $2l$ in \eqref{dpcond} has to be replaced by a
larger integer, determined by the size of the window $f(\omega)$ depends upon. Still,
this gives almost sure absence of eigenvalues for almost every $\theta$.

\medskip

A large number of papers deal with the eigenvalue problem for Schr\"odinger operators
generated by primitive substitutions; for example,
\cite{b4,bbg1,bg1,d3,d4,d5,d6,dp2,hks}.\footnote{There are also papers dealing with
Schr\"odinger operators associated with non-primitive substitutions
\cite{dl5,dol2,dol3,ldo}. The subshifts considered in these papers are, however, linearly
recurrent and hence strictly ergodic, so that the theory is quite similar.} We first
describe general results that can be obtained using the two general methods we discussed
and then turn to some specific examples, where more can be said.

We start with an application of Lemma~\ref{gordon3}. Fix some strictly ergodic subshift
$\Omega$ and define, for $w \in \mathcal{W}_\Omega$, the \textit{index of} $w$ to be
$$
\mathrm{ind}(w) = \sup \{ r \in \Q : w^r \in \mathcal{W}_\Omega \}.
$$
Here, $w^r$ denotes the word $(xy)^m x$, where $m \in \Z^+$, $w = xy$, and $r = m +
|x|/|w|$. The \textit{index of} $\Omega$ is given by
$$
\mathrm{ind}(\Omega) = \sup \{ \mathrm{ind}(w) : w \in \mathcal{W}_\Omega \} \in
[1,\infty].
$$
Then, the following result was shown in \cite{d5} using three-block Gordon.

\begin{theorem}\label{lmpsubst}
Suppose $\Omega$ is generated by a primitive substitution and $\mathrm{ind}(\Omega) > 3$.
Let $f : \Omega \to \R$ be given by $f(\omega) = g(\omega_0)$ with some non-constant
function $g : \mathcal{A} \to \R$. Then, $H_\omega$ has empty point spectrum for
$\mu$-almost every $\omega \in \Omega$.
\end{theorem}

\noindent\textbf{Remarks.} (a) See \cite{d4} for a weaker result, assuming
$\mathrm{ind}(\Omega) \ge 4$.\\
(b) The result extends to the case of a locally constant $f$.\\
(c) Consider the case of the period doubling substitution. Since $s_\mathrm{PD} =
101110101011101110 \ldots$, we see that $\mathrm{ind}(\Omega) \ge \mathrm{ind}(10) \ge
3.5 > 3$. Thus, Theorem~\ref{lmpsubst} implies almost sure absence of eigenvalues,
recovering the main result of \cite{d3}.

\medskip

A substitution belongs to \textit{class P} if there is a palindrome $p$ and, for every $a
\in \mathcal{A}$, a palindrome $q_a$ such that $S(a) = p q_a$. Here, $p$ is allowed to be
the empty word and, if $p$ is not empty, $q_a$ may be the empty word. Clearly, if a
subshift is generated by a class P substitution, it contains arbitrarily long
palindromes. Thus, by Proposition~\ref{hksprop}, it contains uncountably many strongly
palindromic elements. The following result from \cite{hks} is therefore an immediate
consequence.

\begin{theorem}\label{hkssubst}
Suppose $\Omega$ is generated by a primitive substitution $S$ that belongs to class P.
Let $f : \Omega \to \R$ be given by $f(\omega) = g(\omega_0)$ with some non-constant
function $g : \mathcal{A} \to \R$. Then, eigenvalues are generically absent.
\end{theorem}

Notice that the Fibonacci, period doubling, and Thue-Morse subshifts are generated by
class P substitutions. See \cite{hks} for more examples. The Rudin-Shapiro subshift, on
the other hand, is not generated by a class P substitution. In fact, it does not contain
arbitrarily long palindromes \cite{all,b4}.

We mentioned earlier that the proof of Theorem~\ref{hkssubst} cannot give a stronger
result since the set of strongly palindromic sequences is always of zero measure for
substitution subshifts \cite{dz2}. Moreover, it was shown in \cite{d5} that the
three-block Gordon argument cannot prove more than an almost everywhere statement in the
sense that for every minimal aperiodic subshift $\Omega$, there exists an element $\omega
\in \Omega$ such that $\omega$ does not have the infinitely many three block structures
needed for an application of Lemma~\ref{gordon3}. Thus, proofs of uniform results should
use Lemma~\ref{gordon2} in a crucial way. Theorem~\ref{sturmapp} shows that a uniform
result is known in the Fibonacci case, for example, and Lemma~\ref{gordon2} along with
trace map bounds was indeed the key to the proof of this theorem.

Another example for which a uniform result is known is given by the period doubling
substitution \cite{d6}. The trace map bounds are weaker than in the Fibonacci case, but a
combination of two-block and three-block arguments was shown to work. Further
applications of this idea can be found in \cite{ldo}.

The other two examples from Section~\ref{Sec2}, the Thue-Morse and Rudin-Shapiro
substitutions, are not as well understood as Fibonacci and period doubling. Almost sure
or uniform absence of eigenvalues for these cases are open, though expected. Generic
results can be found in \cite{dp2,krri}.

The eigenvalue problem in the context of the other examples mentioned in
Section~\ref{Sec2} has been studied only in a small number of papers. For Arnoux-Rauzy
subshifts, see \cite{dz}; and for interval exchange transformations, see \cite{dog}.

\section{Quantum Dynamics}\label{Sec8}

In this section we focus on the time-dependent Schr\"odinger equation
\begin{equation}\label{tdse}
i \frac{\partial}{\partial t} \psi = H \psi , \quad \psi(0) = \psi_0,
\end{equation}
where $H$ is a Schr\"odinger operator in $\ell^2(\Z)$ with a potential $V : \Z \to \R$,
typically from a strictly ergodic subshift, and $\psi_0 \in \ell^2(\Z)$. By the spectral
theorem, \eqref{tdse} is solved by $\psi(t) = e^{-itH} \psi_0$. Thus, the question we
want to study is the following: Given some potential $V$ and some initial state $\psi_0
\in \ell^2(\Z)$, what can we say about $e^{-itH} \psi_0$ for large times $t$?

\subsection{Spreading of Wavepackets}

Since $\psi_0$ is square-summable, it is in some sense localized near the origin. For
simplicity, one often considers the special case $\psi_0 = \delta_0$---the delta-function
at the origin. With time, $\psi(t)$ will in general spread out in space. Our goal is to
measure this spreading of the wavepacket and relate spreading rates to properties of the
potential. As a general rule of thumb, spreading rates decrease with increased randomness
of the potential. We will make this more explicit below.

A popular way of measuring the spreading of wavepackets is the following. For $p > 0$,
define
\begin{equation}\label{mpo}
\langle |X|_{\psi_0}^p \rangle (T) = \sum_n |n|^p a(n,T),
\end{equation}
where
\begin{equation}\label{ant}
a(n,T)=\frac{2}{T} \int_0^{\infty} e^{-2t/T} | \langle e^{-itH} \psi_0, \delta_n \rangle
|^2 \, dt.
\end{equation}
Clearly, the faster $\langle |X|_{\psi_0}^p \rangle (T)$ grows, the faster $e^{-itH}
\psi_0$ spreads out, at least averaged in time.\footnote{Taking time averages is natural
since the operators of interest in this paper have purely singular continuous spectrum;
compare Wiener's Theorem. While Wiener's Theorem would suggest taking a Ces\`aro time
average, the Abelian time average we choose is more convenient for technical purposes.
The transport exponents are the same for both ways of time averaging.} One typically
wants to obtain power-law bounds on $\langle |X|_{\delta_0}^p \rangle (T)$ and hence it
is natural to define the following quantities: For $p > 0$, define the \textit{upper}
(resp., \textit{lower}) \textit{transport exponent} $\beta^\pm_{\delta_0}(p)$ by
$$
\beta^\pm_{\psi_0}(p)=\operatornamewithlimits{\lim{}^{\text{sup}}_{\text{inf}}}_{T\to\infty}
 \frac{\log \langle |X|_{\psi_0}^p \rangle (T) }{p \, \log T}
$$
Both functions $p \mapsto \beta^\pm_{\psi_0} (p)$ are nondecreasing and obey $0 \le
\beta^-_{\psi_0}(p) \le \beta^+_{\psi_0}(p) \le 1$.

For periodic $V$, $\beta^\pm_{\psi_0} (p) \equiv 1$ (\textit{ballistic transport}); while
for random $V$, $\beta^\pm_{\psi_0} (p) \equiv 0$ (a weak version of \textit{dynamical
localization}---stronger results are known). For $V$'s that are intermediate between
periodic and random, and in particular Sturmian $V$'s, it is expected that the transport
exponents take values between $0$ and $1$.

\subsection{Spectral Measures and Subordinacy Theory}

By the spectral theorem, $\langle e^{-itH} \psi_0 , \psi_0 \rangle = \int e^{-itE} \,
d\mu_{\psi_0}(E)$, where $\mu_{\psi_0}$ is the spectral measure associated with $H$ and
$\psi_0$. Thus, it is natural to investigate quantum dynamical questions by relating them
to properties of the spectral measure corresponding to the initial state. This approach
is classical and the Riemann-Lebesgue Lemma and Wiener's Theorem may be interpreted as
statements in quantum dynamics. The RAGE theorem establishes basic dynamical results in
terms of the standard decomposition of the Hilbert space into pure point, singular
continuous, and absolutely continuous subspaces. We refer the reader to Last's
well-written article \cite{l} for a review of these early results.

The results just mentioned are very satisfactory for initial states whose spectral
measure has an absolutely continuous component. This is, to some extent, also true for
pure point measures. However, if the measure is purely singular continuous, it is
desirable to obtain results that go beyond Wiener's Theorem and the RAGE theorem.

Last also addressed this issue in \cite{l} and proposed a decomposition of spectral
measures with respect to Hausdorff measures. This was motivated by earlier results of
Guarneri \cite{g3} and Combes \cite{c2} who proved dynamical lower bounds for initial
states with uniformly H\"older continuous spectral measures. By approximation with
uniformly H\"older continuous measures, Last proved in \cite{l} that these bounds extend
to measures that are absolutely continuous with respect to a suitable Hausdorff measure:

\begin{theorem}\label{lastthm}
If $\mu_{\psi_0}$ has a non-trivial component that is absolutely continuous with respect
to the $\alpha$-dimensional Hausdorff measure $h^\alpha$ on $\R$, then
\begin{equation}\label{alphalowerbound}
\beta^-_{\psi_0} (p) \ge \alpha \; \text{ for every } p > 0.
\end{equation}
\end{theorem}

\noindent\textbf{Remarks.} (a) Here, $h^\alpha$ is defined by
$$
h^\alpha (S) = \lim_{\delta \to 0} \; \inf_{\text{$\delta$-covers}} \sum |I_m|^\alpha,
$$
where $S \subseteq \R$ is a Borel set and a $\delta$-cover is a countable collection of
intervals $I_m$ of length bounded by $\delta$ such that the union of these intervals
contains the set in question. Note that $h^1$ coincides with Lebesgue measure and $h^0$
is the counting measure.\\
(b) For further developments of quantum dynamical lower bounds in terms of continuity or
dimensionality properties of spectral measures, see \cite{bgt,bt,gsb1,gsb2}.\\
(c) The result and its proof have natural analogues in higher dimensions; see \cite{l}.

\medskip

While a bound like \eqref{alphalowerbound} is nice, it needs to be complemented by
effective methods for verifying the input to Theorem~\ref{lastthm}. In the context of
one-dimensional Schr\"odinger operators, it is always extremely useful to connect a
problem at hand to properties of solutions to the difference equation \eqref{evedet}. The
classical decomposition of spectral measures can be studied via subordinacy theory as
shown by Gilbert and Pearson \cite{gp}; see also \cite{g1,kp}. Subordinacy theory has
proved to be one of the major tools in one-dimensional spectral theory and many important
results have been obtained with its help. Jitomirskaya and Last were able to refine
subordinacy theory to the extent that Hausdorff-dimensional spectral issues can be
investigated in terms of solution behavior \cite{jl0,jl1,jl2}. The key result is the
Jitomirskaya-Last inequality, which explicitly relates the Borel transform of the
spectral measure to solutions in the half-line setting \cite[Theorem~1.1]{jl1}.

Using the maximum modulus principle together with the Jitomirskaya-Last inequality,
Damanik et al.\ then proved the following result for operators on the line \cite{dkl}:

\begin{theorem}\label{dklthm}
Suppose $\Sigma \subseteq \R$ is a bounded set and there are constants $\gamma_1,
\gamma_2$ such that for each $E \in \Sigma$, every solution $u$ of \eqref{evedet} with
$|u(-1)|^2 + |u(0)|^2 = 1$ obeys the estimate
\begin{equation}\label{dklassume}
C_1(E) L^{\gamma_1} \le \left( \sum_{n=1}^L |u(n)|^2 \right)^{1/2} \leq C_2(E)
L^{\gamma_2}
\end{equation}
for $L > 0$ sufficiently large and suitable constants $C_1(E), C_2(E)$. Let $\alpha ={2
\gamma_1}/({\gamma_1 + \gamma_2})$. Then, for any $\psi_0 \in \ell^2(\Z)$, the spectral
measure for the pair $(H,\psi_0)$ is absolutely continuous with respect to $h^\alpha$ on
$\Sigma$. In particular, the bound \eqref{alphalowerbound} holds for every non-trivial
initial state whose spectral measure is supported in $\Sigma$.
\end{theorem}

This shows that suitable bounds for solutions of \eqref{evedet} imply statements on
Hausdorff-dimensional spectral properties, which in turn yield quantum dynamical lower
bounds. There is an extension to multi-dimensional Schr\"odinger operators by Kiselev and
Last \cite{kislas}.

Two remarks are in order. First, while there are some important applications of the
method just presented, proving the required solution estimates is often quite involved.
The number of known applications is therefore still relatively small. Second, dynamical
bounds in terms of Hausdorff-dimensional properties are strictly one-sided. It is not
possible to prove dynamical upper bounds purely in terms of dimensional properties. There
are a number of examples that demonstrate this phenomenon. For example, modifications of
the super-critical almost Mathieu operator lead to spectrally localized operators with
almost ballistic transport \cite{djls,gkt}. Another important example that is spectrally,
but not dynamically, localized is given by the random dimer model \cite{dbg,jss}.

\subsection{Plancherel Theorem}

There is another approach to dynamical bounds that is also based on solution (or rather,
transfer matrix) estimates, which relates dynamics to integrals over Lebesgue measure, as
opposed to integrals over the spectral measure. Compared with the approach discussed
above, it has two main advantages: One can prove both upper and lower bounds in this way,
and the proof of a lower bound is so soft that it applies to a greater number of models.

The key to this approach is a formula due to Kato, which follows quickly from the
Plancherel Theorem:

\begin{lemma}\label{plancherel}
\begin{equation}\label{parsform}
2\pi \int_0^{\infty} e^{-2t/T} | \langle e^{-itH} \psi_0, \delta_n \rangle |^2 \, dt =
\int_{-\infty}^\infty \left|\langle (H - E - \tfrac{i}{T})^{-1} \psi_0, \delta_n \rangle
\right|^2 \, dE.
\end{equation}
\end{lemma}

\begin{proof}
Consider the function
$$
F(t) = \begin{cases} e^{-t/T} \langle e^{-itH} \psi_0, \delta_n \rangle & t \ge 0, \\ 0 &
t < 0. \end{cases}
$$
Using the spectral theorem, it is readily verified that the Fourier transform of $F$
obeys $\widehat{F}(-E) = i \langle (H - E - \tfrac{i}{T})^{-1} \psi_0, \delta_n \rangle$.
Thus, \eqref{parsform} follows if we apply the Plancherel theorem to $F$.
\end{proof}

For simplicity, let us consider the case $\psi_0 = \delta_0$. Note that
$$
u(n) = \langle (H - E - i/T)^{-1} \delta_0, \delta_n \rangle
$$
solves the difference equation \eqref{evedet} (with $E$ replaced by $E + i/T$) away from
the origin and can therefore be studied by means of transfer matrices! In particular, we
may infer a bound from below in terms of $\| A_n^{E + i/T} \|^{-1}$. Thus, upper bounds
on transfer matrix norms are of interest.

\begin{theorem}\label{dtthm}
Suppose that the transfer matrices obey the bound $\| A_n^E \| \le C |n|^{\alpha}$ for
every $n \not= 0$, some fixed energy $E \in \R$ and suitable constants $C,\alpha$. Then,
$$
\beta^-_{\delta_0} (p) \ge \frac{1}{1+2\alpha} - \frac{1+8 \alpha}{p+2\alpha p}
$$
for every $p > 0$.
\end{theorem}

\noindent\textbf{Remarks.} (a) This is the one-energy version of a more general result
due to Damanik and Tcheremchantsev \cite{dt1}. See \cite{dst} for extensions of
\cite{dt1} and supplementary material and \cite{gkt} for related work.\\
(b) An interesting application of Theorem~\ref{dtthm} (and its proof) to the random dimer
model may be found in the paper \cite{jss} by Jitomirskaya et al., which confirms a
prediction of Dunlap et al.\ \cite{dwp}.\\
(c) There is also a version of Theorem~\ref{dtthm} for more general initial states
$\psi_0$ \cite{dls}.\\
(d) The idea of the proof of Theorem~\ref{dtthm} is simple. A Gronwall-type perturbation
argument derives upper bounds on $\| A_n^{\tilde E} \|$ for $\tilde E$ close to $E$ and
$n$ not too large. The right-hand side of \eqref{parsform} may then be estimated from
below by integrating only over a small neighborhood of $E$, where $u$ is controlled by
the upper bound on the transfer matrix. The bound for $\beta^-_{\delta_0} (p)$ then
follows by rather straightforward arguments.\\
(e) The paper \cite{dt2} (using some ideas from \cite{t}) shows that a combination of the
two approaches may sometimes (e.g., in the Fibonacci case) give better bounds.\\
(f) Killip et al.\ used \eqref{parsform} to prove dynamical upper bounds for the slow
part of the wavepacket \cite{kkl}. Their work inspired the use of \eqref{parsform} in
\cite{dt1}.

\medskip

Since \eqref{parsform} is an identity, rather than an inequality, it can be used to bound
$a(n,T)$ from both below and above. Clearly, proving an upper bound is more involved and
will require assumptions that are global in the energy. It was shown by Damanik and
Tcheremchantsev that the following assumption on transfer matrix growth is sufficient to
allow one to infer an upper bound for the transport exponents \cite{dt3}:

\begin{theorem}\label{dtthm2}
Let $K \ge 4$ be such that $\sigma (H) \subseteq [-K+1,K-1]$. Suppose that, for some $C
\in (0,\infty)$ and $\alpha \in (0,1)$, we have
\begin{equation}\label{assumeright}
\int_{-K}^K \left( \max_{1 \le n \le C T^\alpha} \Big\| A^{E+ \tfrac{i}{T}}_n \Big\|^2
\right)^{-1} dE = O(T^{-m})
\end{equation}
and
\begin{equation}\label{assumeleft}
\int_{-K}^K \left( \max_{1 \le -n \le C T^\alpha} \Big\| A^{E+ \tfrac{i}{T}}_n \Big\|^2
\right)^{-1} dE = O(T^{-m})
\end{equation}
for every $m \ge 1$. Then, $\beta^+_{\delta_0} (p) \le \alpha$ for every $p > 0$.
\end{theorem}

\subsection{Applications}

Let us discuss the applications of these general methods to Schr\"odinger operators with
potentials from strictly ergodic subshifts.

We begin with the Fibonacci case. In fact, every approach to quantum dynamical bounds has
been tested on this example and there are many papers proving dynamical results for it;
for example, \cite{d2,d8,dkl,dt1,dt2,dt3,jl2,kkl}.

Upper bounds for transfer matrices were established by Iochum and Testard \cite{it} who
proved, for zero phase, that the norms of the transfer matrices grow no faster than a
power law for every energy in the spectrum. The power can be chosen uniformly on the
spectrum and depends only on the sampling function $f$. Notice that this improves on the
statement that the Lyapunov exponent vanishes on the spectrum. An extension to Sturmian
subshifts whose slope has (essentially) bounded continued fraction coefficients was
obtained by Iochum et al.\ \cite{irt}. Note that upper bounds for transfer matrix norms
yield the input to Theorem~\ref{dtthm} and one half of the input to Theorem~\ref{dklthm}.
The other half of the input to Theorem~\ref{dklthm}, lower bounds for solutions, was
obtained in \cite{d2,jl2}. The proof of these bounds uses the bound for the trace map for
energies from the spectrum, Gordon's two-block lemma, and a mass-reproduction technique
based on cyclic permutations of repeated blocks.\footnote{Using partitions
(cf.~Theorem~\ref{partition}), these solution estimates described in this paragraph can
be shown for all elements of the subshift \cite{dkl,dl2}.}

\begin{theorem}\label{fibdynlowerthm}
Let $V(n) = \lambda \chi_{[1-\theta,1)}(n \theta)$, where $\theta = (\sqrt{5}-1)/2$ and
$\lambda > 0$. Then,
\begin{equation}\label{dt2lower}
\beta^-_{\delta_0} (p) \ge \begin{cases} \frac{p+2\kappa}{(p+1)(\alpha + \kappa + 1/2)} & p \le 2 \alpha+1, \\
\frac{1}{\alpha+1} & p>2 \alpha +1. \end{cases},
\end{equation}
where $\kappa$ is an absolute constant {\rm (}$\kappa \approx 0.0126${\rm )} and $\alpha
\asymp \log \lambda$.
\end{theorem}

\noindent\textbf{Remarks.} (a) In the form stated, the result is from \cite{dt2}. The
bound \eqref{dt2lower} is the best known dynamical lower bound for the Fibonacci operator
and is a culmination of the sequence of works \cite{d2,dkl,dt1,jl2,kkl} leading up to
\cite{dt2}.\\
(b) When we write $\alpha \asymp \log \lambda$, we mean that $\alpha$ is a positive
$\lambda$-dependent quantity that satisfies $C_1 \log \lambda \le \alpha \le C_2 \log
\lambda$ for positive constants $C_1,C_2$ and all large $\lambda$. See \cite{dt2} for the
explicit dependence of $\alpha$ on $\lambda$.

\medskip

To apply Theorem~\ref{dtthm2} to the Fibonacci operator, one has to prove the estimates
\eqref{assumeright} and \eqref{assumeleft}. This was done in \cite{dt3}. Let us describe
the main idea. Clearly, to prove the desired lower bounds for transfer matrix norms, it
suffices to prove lower bounds for transfer matrix traces. We know a way to establish
such lower bounds: Lemma~\ref{fibdicho}. Since all relevant energies in
\eqref{assumeright} and \eqref{assumeleft} are non-real, we know that the trace map will
eventually enter the escape region described in Lemma~\ref{fibdicho}. The point is to
control the number of iterates it takes for this to occur. To this end, define the
complex analogue of the set $\sigma_k$ from Section~\ref{Sec6} by
$$
\sigma^\C_k = \{ z \in \C : |x_k(z)| \le 1 \}.
$$
Notice that the $x_k$'s are polynomials and hence defined for all complex $z$. As before,
being in the complement of two consecutive $\sigma^\C_k$'s is a sufficient condition for
escape at an explicit rate; compare Lemma~\ref{fibdicho}, whose proof extends to complex
energies. It is therefore useful to bound the imaginary width of these sets from above.
This will give an upper bound on the number of iterates it takes at a given energy to
enter the escape region. For $\lambda$ sufficiently large, the connected components of
$\sigma^\C_k$ can be studied with the help of Koebe's Distortion Theorem; see \cite{dt3}
for details. The resulting dynamical upper bound has the same asymptotics for large
$\lambda$ as the lower bound above:

\begin{theorem}\label{fibdynupperthm}
Let $V(n) = \lambda \chi_{[1-\theta,1)}(n \theta)$, where $\theta = (\sqrt{5}-1)/2$ and
$\lambda \ge 8$. Then,
$$
\beta^+_{\delta_0} (p) \le \tilde \alpha \quad \text{ for every } p > 0,
$$
where $\tilde \alpha \in (0,1)$ and $\tilde \alpha \asymp (\log \lambda)^{-1}$.
\end{theorem}

In particular, for the Fibonacci operator with $\lambda \ge 8$, all transport exponents
$\{ \beta^\pm_{\delta_0} (p) \}_{p > 0}$ are strictly between zero and one.

The dynamical lower bounds have been established for more general models; see
\cite{d2,dkl,dl4,dls,dst,dt1}. On the other hand, Theorem~\ref{fibdynupperthm} is the
only explicit result of this kind, but as mentioned in \cite{dt3}, the ideas of \cite{d8}
should permit one to extend this theorem to more general slopes and all elements of the
subshift.

\section{CMV Matrices Associated with Subshifts}

Given a strictly ergodic subshift $\Omega$ and a continuous/locally constant function $f
: \Omega \to \D$, we can define $\alpha_n(\omega) = f(T^n \omega)$ for $n \in \Z$ and
$\omega \in \Omega$. Let $\mathcal{C}_\omega$ be the CMV matrix associated with
Verblunsky coefficients $\{\alpha_n(\omega)\}_{n \ge 0}$ and $\mathcal{E}_\omega$ the
extended CMV matrix associated with Verblunsky coefficients $\{\alpha_n(\omega)\}_{n \in
\Z}$. That is, with $\rho_n(\omega) = (1 - |\alpha_n(\omega)|)^{-1/2}$,
$\mathcal{C}_\omega$ is given by
$$
\begin{pmatrix}
{}& \bar\alpha_0(\omega) & \bar\alpha_1(\omega) \rho_0(\omega) & \rho_1(\omega) \rho_0(\omega)
& 0 & 0 & \dots & {} \\
{}& \rho_0(\omega) & -\bar\alpha_1(\omega) \alpha_0(\omega) & -\rho_1(\omega) \alpha_0(\omega)
& 0 & 0 & \dots & {} \\
{}& 0 & \bar\alpha_2(\omega) \rho_1(\omega) & -\bar\alpha_2(\omega) \alpha_1(\omega) &
\bar\alpha_3(\omega) \rho_2(\omega) & \rho_3(\omega) \rho_2(\omega) & \dots & {} \\
{}& 0 & \rho_2(\omega) \rho_1(\omega) & -\rho_2(\omega) \alpha_1(\omega) & -\bar\alpha_3(\omega)
\alpha_2(\omega) & -\rho_3(\omega) \alpha_2(\omega) & \dots & {} \\
{}& 0 & 0 & 0 & \bar\alpha_4(\omega) \rho_3(\omega) & -\bar\alpha_4(\omega) \alpha_3(\omega)
& \dots & {} \\
{}& \dots & \dots & \dots & \dots & \dots & \dots & {}
\end{pmatrix}
$$
and $\mathcal{E}_\omega$ is the analogous two-sided infinite matrix. See \cite{simon,
simon2} for more information on CMV and extended CMV matrices.

For these unitary operators in $\ell^2$, we can ask questions similar to the ones
considered above in the context of Schr\"odinger operators. That is, is the spectrum of
zero Lebesgue measure, are spectral measures purely singular continuous, etc. Since we
are dealing with ergodic models, it is more natural to consider the whole-line situation.
On the other hand, from the point of view of orthogonal polynomials on the unit circle,
the half-line situation is more relevant. The zero-measure property is independent of the
setting, whereas the spectral type for half-line models is almost always (i.e., when the
``boundary condition'' is varied) pure point as soon as zero-measure spectrum is
established. The latter statement follows quickly from spectral averaging; compare
\cite[Theorem~10.2.2]{simon2}. Thus, the key problem for CMV matrices associated with
subshifts is proving zero-measure spectrum. In fact, Simon conjectured the following; see
\cite[Conjecture~12.8.2]{simon2}.

\medskip

\noindent\textbf{Simon's Subshift Conjecture.} Suppose $\mathcal{A}$ is a subset of $\D$,
the subshift $\Omega$ is minimal and aperiodic and let  $f : \Omega \to \D$, $f(\omega)=
\omega(0)$. Then, $\Sigma$ has zero Lebesgue measure.

\medskip

Here, $\Sigma$ is the common spectrum of the operators $\mathcal{E}_\omega$, $\omega \in
\Omega$. Equivalently, it is the common essential spectrum of $\mathcal{C}_\omega$,
$\omega \in \Omega$.

Simon proved the zero-measure property for the Fibonacci case by means of the trace map
approach \cite[Section~12.8]{simon2}. Since the approach based on the Boshernitzan
condition has a wider scope in the Schr\"odinger case, it is natural to try and extend it
to the CMV case. This was done by Damanik and Lenz in \cite{dl12} where the following
result was shown.

\begin{theorem}\label{dlvcthm}
Suppose the subshift $\Omega$ is aperiodic and satisfies the Boshernitzan condition. Let
$f : \Omega \to \D$ be locally constant. Then, $\Sigma$ has zero Lebesgue measure.
\end{theorem}

This proves Simon's Subshift Conjecture for a large number of models since we saw above
that many of the prominent aperiodic subshifts satisfy the Boshernitzan condition.

Regarding the spectral type, it should not be hard to extend the material from
Section~\ref{Sec7} to the CMV case. This will imply purely singular continuous spectrum
for $\mathcal{E}_\omega$ for many subshifts $\Omega$ and many (generic, almost all, all)
$\omega \in \Omega$. However, as was noted above, the Aleksandrov measures associated
with $\mathcal{C}_\omega$ will almost surely be pure point whenever Theorem~\ref{dlvcthm}
applies.

Quantum dynamics, on the other hand, is less natural in the CMV case than in the
Schr\"odinger case, and has not really been studied.\footnote{Simon did extend the
Jitomirskaya-Last theory to OPUC in \cite{simon2}. This theory has its roots in quantum
dynamics; compare \cite{djls,g3,jl1,jl2,l}.} Most of the ideas leading to the results
presented in Section~\ref{Sec8} should have CMV counterparts. In particular, it should be
possible to prove absolute continuity of spectral measures with respect to suitable
Hausdorff measures for extended CMV matrices over Fibonacci-like subshifts.

\section{Concluding Remarks}

The material presented in this survey is motivated by and closely related to the theory
of quasicrystals; compare, for example, \cite{baamoo,moo}. More specifically, the surveys
\cite{d7,s7} deal with the Fibonacci operator and its generalizations and the interested
reader may find references to the original physics literature in those papers.

Regarding future research in this field, it would be interesting to see how far one can
take the philosophy that potentials taking finitely many values preclude localization
phenomena. Since the Bernoulli Anderson model is localized \cite{ckm}, this cannot hold
in full generality. On the other hand, Gordon potentials are much more prevalent in the
subshift case than in the uniformly almost periodic case. Moreover, for smooth
quasi-periodic potentials, it is expected that the Lyapunov exponent is positive at all
energies if the coupling is large enough. This is known for trigonometric potentials
\cite{herman}, analytic potentials \cite{bgold,golsch,sorspe}, and Gevrey potentials
\cite{klein}. See also \cite{bjerk,chan} for recent results in the $C^r$ category. These
potentials should be contrasted with those coming from quasi-periodic subshifts
satisfying the Boshernitzan condition. The Boshernitzan condition is independent of the
coupling constant and yields vanishing Lyapunov exponent throughout the spectrum. Since
it is satisfied on a dense set of sampling (step-)functions, upper-semicontinuity
arguments allow one to derive surprising phenomena that hold generically in the $C^0$
category \cite{bdj}.

To shed some light on this, it should be helpful to analyze more examples. That is, take
one of the popular base transformations of the torus (e.g., shifts, skew-shifts, or
expanding maps) and define an ergodic family of potentials by choosing a sampling
function on the torus that takes finitely many values. These models, with the exception
of rotations of the circle, are not well understood! There is a serious issue about the
competition between the flat pieces of the sampling function and the randomness
properties of the base transformation (expressed, e.g., in terms of mixing properties).
For example, take a $1$-periodic step function $f$ and consider $V_\omega (n) = \lambda
f(2^n \omega)$, $\lambda > 0$, $\omega \in [0,1)$. Is it true that the Lyapunov exponent
is positive? For all $\lambda$'s or all large $\lambda$'s? For all energies or all but
finitely many?


\begin{thebibliography}{9999}

\bibitem{ab2} P.\ Alessandri and V.\ Berth\'e, Three distance theorems and combinatorics on words,
\textit{Enseign.\ Math.} \textbf{44} (1998), 103--132

\bibitem{all} J.-P.\ Allouche, Schr\"odinger operators with Rudin-Shapiro potentials are
not palindromic, \textit{J.\ Math.\ Phys.} \textbf{38} (1997), 1843--1848

\bibitem{ap} J.-P.\ Allouche and J.\ Peyri\`{e}re, Sur une formule de r\'{e}currence sur les
traces de produits de matrices associ\'{e}s a certaines substitutions, \textit{C.\ R.\
Acad.\ Sci.\ Paris} {\bf 302} (1986), 1135--1136

\bibitem{ab} Y.\ Avishai and D.\ Berend, Trace maps for arbitrary substitution sequences, \textit{J.\
Phys.\ A} {\bf 26} (1993), 2437--2443

\bibitem{abg} Y.\ Avishai, D.\ Berend, and D.\ Glaubman, Minimum-dimension trace maps for
substitution sequences, \textit{Phys.\ Rev.\ Lett.} {\bf 72} (1994), 1842--1845

\bibitem{avsim} J.\ Avron and B.\ Simon, Singular continuous spectrum for a class of almost periodic
Jacobi matrices, \textit{Bull.\ Amer.\ Math.\ Soc.} \textbf{6} (1982), 81--85

\bibitem{b4} M.\ Baake, A note on palindromicity, \textit{Lett.\ Math.\ Phys.} {\bf 49} (1999),
217--227

\bibitem{bgj} M.\ Baake, U.\ Grimm, and D.\ Joseph, Trace maps, invariants, and some of their
applications, \textit{Int.\ J.\ Mod.\ Phys.\ B} {\bf 7} (1993), 1527--1550

\bibitem{baamoo} M.\ Baake and R.\ Moody (Editors), \textit{Directions in Mathematical
Quasicrystals}, American Mathematical Society, Providence (2000)

\bibitem{br} M.\ Baake and J.\ Roberts, Reversing symmetry group of $\mathrm{GL}(2,\Z)$ and
$\mathrm{PGL}(2,\Z)$ matrices with connections to cat maps and trace maps, \textit{J.\
Phys.\ A} {\bf 30} (1997), 1549--1573

\bibitem{bgt} J.-M.\ Barbaroux, F.\ Germinet, and S.\ Tcheremchantsev, Fractal dimensions and
the phenomenon of intermittency in quantum dynamics, \textit{Duke Math.\ J.} \textbf{110}
(2001), 161--193

\bibitem{bt} J.-M.\ Barbaroux and S.\ Tcheremchantsev, Universal lower bounds for quantum diffusion,
\textit{J.\ Funct.\ Anal.} \textbf{168} (1999), 327--354

\bibitem{b2} J.\ Bellissard, Spectral properties of Schr\"odinger's operator with a Thue-Morse
potential, in \textit{Number Theory and Physics} (\textit{Les Houches, 1989}), Springer,
Berlin (1990), 140--150

\bibitem{bbg1} J.\ Bellissard, A.\ Bovier, and J.-M.\ Ghez, Spectral properties of a tight binding
Hamiltonian with period doubling potential, \textit{Commun.\ Math.\ Phys.} {\bf 135}
(1991), 379--399

\bibitem{bist} J.\ Bellissard, B.\ Iochum, E.\ Scoppola, and D.\ Testard, Spectral properties of
one-dimensional quasicrystals, \textit{Commun.\ Math.\ Phys.} {\bf 125} (1989), 527--543

\bibitem{berstel} J.\ Berstel, Recent results in Sturmian words, in \textit{Developments in Language
Theory}, World Scientific, Singapore (1996), 13--24

\bibitem{bjerk} K.\ Bjerkl\"ov, Positive Lyapunov exponent and minimality for a class of 1-d
quasi-periodic Schr\"odinger equations, to appear in \textit{Ergod.\ Th.\ \& Dynam.\
Sys.}

\bibitem{bdj} K.\ Bjerkl\"ov, D.\ Damanik, and R.\ Johnson, Lyapunov exponents of continuous
Schr\"odinger cocycles over irrational rotations, Preprint (2005)

\bibitem{Bosh1} M.\ Boshernitzan, A condition for minimal interval exchange
maps to be uniquely ergodic, \textit{Duke Math.\ J.} {\bf  52} (1985), 723--752

\bibitem{Bosh2} M.\ Boshernitzan, A condition for unique ergodicity of minimal
symbolic flows, \textit{Ergod.\ Th.\ \& Dynam.\ Sys.} {\bf 12} (1992), 425--428

\bibitem{bgold} J.\ Bourgain and M.\ Goldstein, On nonperturbative localization with quasi-periodic
potential, \textit{Ann.\ of Math.} {\bf 152} (2000), 835--879

\bibitem{bg1} A.\ Bovier and J.-M.\ Ghez, Spectral properties of one-dimensional Schr\"odinger
operators with potentials generated by substitutions, \textit{Commun.\ Math.\ Phys.} {\bf
158} (1993), 45--66; Erratum: \textit{Commun.\ Math.\ Phys.} {\bf 166} (1994), 431--432

\bibitem{ckm} R.\ Carmona, A.\ Klein, and F.\ Martinelli, Anderson localization for Bernoulli and
other singular potentials, \textit{Commun.\ Math.\ Phys.} {\bf 108} (1987), 41--66

\bibitem{cl} R.\ Carmona and J.\ Lacroix, \textit{Spectral Theory of Random Schr\"odinger Operators},
Birkh\"auser, Boston (1990)

\bibitem{c} M.\ Casdagli, Symbolic dynamics for the renormalization group of a quasiperiodic
Schr\"odinger equation, \textit{Commun.\ Math.\ Phys.} {\bf 107} (1986), 295--318

\bibitem{chan} J.\ Chan, Method of variations of potential of quasi-periodic
Schr\"odinger equation, Preprint (2005)

\bibitem{chna} J.\ Choksi and M.\ Nadkarni, Genericity of certain classes of unitary and self-adjoint
operators, \textit{Canad.\ Math.\ Bull.} \textbf{41} (1998), 137--139

\bibitem{c2} J.-M.\ Combes, Connections between quantum dynamics and spectral properties of
time-evolution operators, in \textit{Differential Equations with Applications to
Mathematical Physics}, Academic Press, Boston (1993), 59--68

\bibitem{d2} D.\ Damanik, $\alpha$-continuity properties of one-dimensional quasicrystals,
\textit{Commun.\ Math.\ Phys.} {\bf 192} (1998), 169--182

\bibitem{d3} D.\ Damanik, Singular continuous spectrum for the period doubling Hamiltonian on a
set of full measure, \textit{Commun.\ Math.\ Phys.} {\bf 196} (1998), 477--483

\bibitem{d4} D.\ Damanik, Singular continuous spectrum for a class of substitution
Hamiltonians, \textit{Lett.\ Math.\ Phys.} {\bf 46} (1998), 303--311

\bibitem{d5} D.\ Damanik, Singular continuous spectrum for a class of substitution Hamiltonians
II., \textit{Lett.\ Math.\ Phys.} \textbf{54} (2000), 25--31

\bibitem{d7} D.\ Damanik, Gordon-type arguments in the spectral theory of
one-dimensional quasicrystals, in \textit{Directions in Mathematical Quasicrystals},
American Mathematical Society, Providence (2000), 277--305

\bibitem{d6} D.\ Damanik, Uniform singular continuous spectrum for the period doubling Hamiltonian,
\textit{Ann.\ Henri Poincar\'e} \textbf{2} (2001), 101--108

\bibitem{d8} D.\ Damanik, Dynamical upper bounds for one-dimensional quasicrystals, \textit{J.\ Math.\ Anal.\
Appl.} \textbf{303} (2005), 327--341

\bibitem{dgr} D.\ Damanik, J.-M.\ Ghez, and L.\ Raymond, A palindromic half-line criterion for absence
of eigenvalues and applications to substitution Hamiltonians, \textit{Ann.\ Henri
Poincar\'{e}} \textbf{2} (2001), 927--939

\bibitem{dkl} D.\ Damanik, R.\ Killip, and D.\ Lenz, Uniform spectral properties of one-dimensional
quasicrystals. III.~$\alpha$-continuity, \textit{Commun.\ Math.\ Phys.} {\bf 212} (2000),
191--204

\bibitem{dl1} D.\ Damanik and D.\ Lenz, Uniform spectral properties of one-dimensional
quasicrystals, I.~Absence of eigenvalues, \textit{Commun.\ Math.\ Phys.} {\bf 207}
(1999), 687--696

\bibitem{dl2} D.\ Damanik and D.\ Lenz, Uniform spectral properties of one-dimensional
quasicrystals, II.~The Lyapunov exponent, \textit{Lett.\ Math.\ Phys.} {\bf 50} (1999),
245--257

\bibitem{dl11} D.\ Damanik and D.\ Lenz, Linear repetitivity, I.\ Uniform subadditive ergodic
theorems and applications, \textit{Discrete Comput.\ Geom.} {\bf 26} (2001), 411--428

\bibitem{dl10} D.\ Damanik and D.\ Lenz, The index of Sturmian sequences, \textit{European
J.\ Combin.} \textbf{23} (2002), 23--29

\bibitem{dl9} D.\ Damanik and D.\ Lenz, Powers in Sturmian sequences, \textit{European J.\
Combin.} \textbf{24} (2003), 377--390

\bibitem{dl4} D.\ Damanik and D.\ Lenz, Uniform spectral properties of one-dimensional
quasicrystals, IV.~Quasi-Sturmian potentials, \textit{J.\ Anal.\ Math.} {\bf 90} (2003),
115--139

\bibitem{dl13} D.\ Damanik and D.\ Lenz, Half-line eigenfunction estimates and purely singular
continuous spectrum of zero Lebesgue measure, \textit{Forum Math.} \textbf{16} (2004),
109--128

\bibitem{dl5} D.\ Damanik and D.\ Lenz, Substitution dynamical systems: Characterization
of linear repetitivity and applications, to appear in \textit{J.\ Math.\ Anal.\ Appl.}

\bibitem{dl7}  D.\ Damanik and D.\ Lenz, A condition of Boshernitzan and uniform convergence
in the Multiplicative Ergodic Theorem, to appear in \textit{Duke Math.\ J.}

\bibitem{dl12}  D.\ Damanik and D.\ Lenz, Uniform Szeg\H{o} cocycles over strictly ergodic
subshifts, Preprint (2005)

\bibitem{dl8}  D.\ Damanik and D.\ Lenz, Zero-measure Cantor spectrum for Schr\"odinger
operators with low-complexity potentials, Preprint (2005)

\bibitem{dls} D.\ Damanik, D.\ Lenz, and G.\ Stolz, Lower transport bounds for one-dimensional continuum
Schr\"odinger operators, Preprint (2004), arXiv/math-ph/0410062

\bibitem{dst} D.\ Damanik, A.\ S\"ut\H{o}, and S.\ Tcheremchantsev, Power-law bounds on transfer
matrices and quantum dynamics in one dimension II, \textit{J.\ Funct.\ Anal.}
\textbf{216} (2004), 362--387

\bibitem{dt1} D.\ Damanik and S.\ Tcheremchantsev, Power-law bounds on transfer matrices and
quantum dynamics in one dimension, \textit{Commun.\ Math.\ Phys.} {\bf 236} (2003),
513--534

\bibitem{dt2} D.\ Damanik and S.\ Tcheremchantsev, Scaling estimates for solutions and dynamical lower
bounds on wavepacket spreading, to appear in \textit{J.\ Anal.\ Math.}

\bibitem{dt3} D.\ Damanik and S.\ Tcheremchantsev, Upper bounds in quantum dynamics,
Preprint (2005), arXiv/math-ph/0502044

\bibitem{dz} D.\ Damanik and L.\ Q.\ Zamboni, Combinatorial properties of Arnoux-Rauzy
subshifts and applications to Schr\"odinger operators, \textit{Rev.\ Math.\ Phys.}
\textbf{15} (2003), 745--763

\bibitem{dz2} D.\ Damanik and D.\ Zare, Palindrome complexity bounds for primitive substitution
sequences, \textit{Discrete Math.} \textbf{222} (2000), 259--267

\bibitem{dbg} S.\ De Bi\`evre and F.\ Germinet, Dynamical localization for the random dimer Schr\"odinger
operator, \textit{J.\ Stat.\ Phys.} \textbf{98} (2000), 1135--1148

\bibitem{djls} R.\ del Rio, S.\ Jitomirskaya, Y.\ Last, and B.\ Simon, Operators with singular continuous
spectrum. IV.~Hausdorff dimensions, rank one perturbations, and localization, \textit{J.\
Anal.\ Math.} \textbf{69} (1996), 153--200

\bibitem{dp1} F.\ Delyon and D.\ Petritis, Absence of localization in a class of Schr\"odinger
operators with quasiperiodic potential, \textit{Commun.\ Math.\ Phys.} {\bf 103} (1986),
441--444

\bibitem{dp2} F.\ Delyon and J.\ Peyri\`{e}re, Recurrence of the eigenstates of a Schr\"odinger
operator with automatic potential, \textit{J.\ Stat.\ Phys.} {\bf 64} (1991), 363--368

\bibitem{dog} C.\ de Oliveira and C.\ Gutierrez, Almost periodic Schr\"odinger operators along
interval exchange transformations, \textit{J.\ Math.\ Anal.\ Appl.} \textbf{283} (2003),
570--581

\bibitem{dol2} C.\ de Oliveira and M.\ Lima, A nonprimitive substitution Schr\"odinger
operator with generic singular continuous spectrum, \textit{Rep.\ Math.\ Phys.}
\textbf{45} (2000), 431--436

\bibitem{dol3} C.\ de Oliveira and M.\ Lima, Singular continuous spectrum for a class of
nonprimitive substitution Schr\"odinger operators, \textit{Proc.\ Amer.\ Math.\ Soc.}
\textbf{130} (2002), 145--156

\bibitem{djp} X.\ Droubay, J.\ Justin, and G.\ Pirillo, Epi-Sturmian words and
some constructions of de Luca and Rauzy, \textit{Theoret.\ Comput.\ Sci.} {\bf 255}
(2001), 539--553

\bibitem{dwp} D.\ Dunlap, H.-L.\ Wu, and P.\ Phillips, Absence of localization
in a random-dimer model, \textit{Phys.\ Rev.\ Lett.} \textbf{65} (1990), 88--91

\bibitem{dur} F.\ Durand, Linearly recurrent subshifts have a finite number of non-periodic
subshift factors, \textit{Ergod.\ Th.\ \& Dynam.\ Sys.} \textbf{20} (2000), 1061--1078

\bibitem{duhosk} F.\ Durand, B.\ Host, and C.\ Skau, Substitutional dynamical systems, Bratteli
diagrams and dimension groups, \textit{Ergod.\ Th.\ \& Dynam.\ Sys.} \textbf{19} (1999),
953--993

\bibitem{gkt} F.\ Germinet, A.\ Kiselev, and S.\ Tcheremchantsev, Transfer matrices and transport for
Schr\"odinger operators, \textit{Ann.\ Inst.\ Fourier} \textbf{54} (2004), 787--830

\bibitem{g1} D.\ Gilbert, On subordinacy and analysis of the spectrum of Schr\"odinger
operators with two singular endpoints, \textit{Proc.\ Roy.\ Soc.\ Edinburgh A} {\bf 112}
(1989), 213--229

\bibitem{gp} D.\ Gilbert and D.\ Pearson, On subordinacy and analysis of the spectrum of
one-dimensional Schr\"odinger operators, \textit{J.\ Math.\ Anal.\ Appl.} {\bf 128}
(1987), 30--56

\bibitem{golsch} M.\ Goldstein and W.\ Schlag, H\"older continuity of the integrated density of states
for quasi-periodic Schr\"odinger equations and averages of shifts of subharmonic
functions, \textit{Ann.\ of Math.} {\bf 154} (2001), 155--203

\bibitem{g2} A.\ Gordon, On the point spectrum of the one-dimensional Schr\"odinger operator,
\textit{Usp.\ Math.\ Nauk.} {\bf 31} (1976), 257--258

\bibitem{g3} I.\ Guarneri, Spectral properties of quantum diffusion on discrete lattices,
\textit{Europhys.\ Lett.} {\bf 10} (1989), 95--100

\bibitem{gsb1} I.\ Guarneri and H.\ Schulz-Baldes, Lower bounds on wave packet propagation by
packing dimensions of spectral measures, \textit{Math.\ Phys.\ Electron.\ J.} {\bf 5}
(1999), paper~1

\bibitem{gsb2} I.\ Guarneri and H.\ Schulz-Baldes, Intermittent lower bound on quantum
diffusion, \textit{Lett.\ Math.\ Phys.} \textbf{49} (1999), 317--324

%\bibitem{ga} G.\ Gumbs and M.\ Ali, Dynamical maps, Cantor spectra, and localization for
%Fibonacci and related quasiperiodic lattices, \textit{Phys.\ Rev.\
%Lett.} {\bf 60} (1988), 1081--1084

\bibitem{hk} F.\ Hahn and Y.\ Katznelson, On the entropy of uniquely ergodic transformations,
\textit{Trans.\ Amer.\ Math.\ Soc.} \textbf{126} (1967), 335--360

\bibitem{herman} M.\ Herman, Une m\'{e}thode pour minorer les exposants de Lyapunov
et quelques exemples montrant the caract\`{e}re local d'un th\'{e}or\`{e}me d'Arnold et
de Moser sur le tore de dimension $2$, \textit{Comment.\ Math.\ Helv} {\bf 58} (1983),
4453--502

\bibitem{h} A.\ Hof, Some remarks on discrete aperiodic Schr\"odinger operators, \textit{J.\ Stat.\
Phys.} {\bf 72} (1993), 1353--1374

\bibitem{hks} A.\ Hof, O.\ Knill, and B.\ Simon, Singular continuous spectrum for palindromic
Schr\"odinger operators, \textit{Commun.\ Math.\ Phys.} {\bf 174} (1995), 149--159

\bibitem{irt} B.\ Iochum, L.\ Raymond, and D.\ Testard, Resistance of one-dimensional
quasicrystals, \textit{Physica A} {\bf 187} (1992), 353--368

\bibitem{it} B.\ Iochum and D.\ Testard, Power law growth for the resistance in the Fibonacci
model, \textit{J.\ Stat.\ Phys.} {\bf 65} (1991), 715--723

\bibitem{jl0} S.\ Jitomirskaya and Y.\ Last, Dimensional Hausdorff properties of singular
continuous spectra, \textit{Phys.\ Rev.\ Lett.} {\bf 76} (1996), 1765--1769

\bibitem{jl1} S.\ Jitomirskaya and Y.\ Last, Power-law subordinacy and singular spectra.
I. Half-line operators, \textit{Acta Math.} {\bf 183} (1999), 171--189

\bibitem{jl2} S.\ Jitomirskaya and Y.\ Last, Power-law subordinacy and singular spectra.
II. Line operators, \textit{Commun.\ Math.\ Phys.} {\bf 211} (2000), 643--658

\bibitem{jss} S.\ Jitomirskaya, H.\ Schulz-Baldes, and G.\ Stolz, Delocalization in random
polymer models, \textit{Commun.\ Math.\ Phys.} {\bf 233} (2003), 27--48

\bibitem{js} S.\ Jitomirskaya and B.\ Simon, Operators with singular continuous spectrum:
III.~Almost periodic Schr\"odinger operators, \textit{Commun.\ Math.\ Phys.} {\bf 165}
(1994), 201--205

\bibitem{john} R.\ Johnson, Exponential dichotomy, rotation number, and linear differential
operators with bounded coefficients, \textit{J.\ Differential Equations} \textbf{61}
(1986), 54--78

\bibitem{juspir} J.\ Justin and G.\ Pirillo, Fractional powers in Sturmian words, \textit{Theoret.\
Comput.\ Sci.} \textbf{255} (2001), 363--376

\bibitem{juspir2} J.\ Justin and G.\ Pirillo, Episturmian words and episturmian morphisms,
\textit{Theoret.\ Comput.\ Sci.} {\bf 276} (2002), 281--313

\bibitem{k1} M.\ Kaminaga, Absence of point spectrum for a class of discrete Schr\"odinger
operators with quasiperiodic potential, \textit{Forum Math.} {\bf 8} (1996), 63--69

\bibitem{Ke1} M.\ Keane, Interval exchange transformations, \textit{Math.\ Z.} {\bf
141} (1975), 25--31

\bibitem{Ke2} M.\ Keane, Non-ergodic interval exchange transformations, \textit{Israel
J.\ Math.} {\bf 26} (1977), 188--196

\bibitem{KN} H.\ B.\ Keynes and D.\ Newton, A minimal, non-uniquely ergodic interval exchange
transformation, \textit{Math.\ Z.} {\bf 148} (1976), 101--105

\bibitem{kp} S.\ Khan and D.\ Pearson, Subordinacy and spectral theory for infinite
matrices, \textit{Helv.\ Phys.\ Acta} {\bf 65} (1992), 505--527

\bibitem{khin} A.\ Khintchine, \textit{Continued Fractions}, Dover, Mineola (1997)

\bibitem{kkl} R.\ Killip, A.\ Kiselev, and Y.\ Last, Dynamical upper bounds on wavepacket
spreading, \textit{Amer.\ J.\ Math.} {\bf 125} (2003), 1165--1198

\bibitem{kislas} A.\ Kiselev and Y.\ Last, Solutions, spectrum, and dynamics for Schr\"odinger
operators on infinite domains, \textit{Duke Math.\ J.} \textbf{102} (2000), 125--150

\bibitem{klein} S.\ Klein, Anderson localization for the discrete one-dimensional quasi-periodic
Schr\"odinger operator with potential defined by a Gevrey-class function, \textit{J.\
Funct.\ Anal.} \textbf{218} (2005), 255--292

\bibitem{kkt} M.\ Kohmoto, L.\ Kadanoff, and C.\ Tang, Localization problem in one dimension:
Mapping and escape, \textit{Phys.\ Rev.\ Lett.} {\bf 50} (1983), 1870--1872

\bibitem{kn} M.\ Kol\'{a}r and F.\ Nori, Trace maps of general substitutional sequences,
\textit{Phys.\ Rev.\ B} {\bf 42} (1990), 1062--1065

\bibitem{k2} S.\ Kotani, Ljapunov indices determine absolutely continuous spectra of stationary
random one-dimensional Schr\"odinger operators, in \textit{Stochastic Analysis}
(\textit{Katata/Kyoto, 1982}), North Holland, Amsterdam (1984), 225--247

\bibitem{k3} S.\ Kotani, Jacobi matrices with random potentials taking finitely many values,
\textit{Rev.\ Math.\ Phys.} {\bf 1} (1989), 129--133

\bibitem{k4} S.\ Kotani, Generalized Floquet theory for stationary Schr\"odinger operators in one
dimension, \textit{Chaos Solitons Fractals} \textbf{8} (1997), 1817--1854

\bibitem{krri} L.\ Kroon and R.\ Riklund, Absence of localization in a model with
correlation measure as a random lattice, \textit{Phys.\ Rev.\ B} \textbf{69} (2004),
paper~094204 (5 pages)

\bibitem{lagple} J.\ Lagarias and P.\ Pleasants, Repetitive Delone sets and quasicrystals,
\textit{Ergod.\ Th.\ \& Dynam.\ Sys.} \textbf{23} (2003), 831--867

\bibitem{l} Y.\ Last, Quantum dynamics and decompositions of singular continuous spectra, \textit{J.\
Funct.\ Anal.} {\bf 142} (1996), 406--445

\bibitem{ls} Y.\ Last and B.\ Simon, Eigenfunctions, transfer matrices, and absolutely
continuous spectrum of one-dimensional Schr\"odinger operators, \textit{Invent.\ Math.}
{\bf 135} (1999), 329--367

\bibitem{len1} D.\ Lenz, Uniform ergodic theorems on subshifts over a finite alphabet,
\textit{Ergod.\ Th.\ \& Dynam.\ Sys.} {\bf 22} (2002), 245--255

\bibitem{len2} D.\ Lenz, Singular continuous spectrum of Lebesgue measure zero for
one-dimensional quasicrystals, \textit{Commun.\ Math.\ Phys.} {\bf 227} (2002), 119--130

\bibitem{lest} D.\ Lenz and P.\ Stollmann, Generic sets in spaces of measures and generic singular
continuous spectrum for Delone Hamiltonians, to appear in \textit{Duke Math.\ J.}

\bibitem{ldo} M.\ Lima and C.\ de Oliveira, Uniform Cantor singular continuous spectrum for
nonprimitive Schr\"odinger operators, \textit{J.\ Statist.\ Phys.} \textbf{112} (2003),
357--374

\bibitem{ltww} Q.-H.\ Liu, B.\ Tan, Z.-X.\ Wen, and J.\ Wu, Measure zero spectrum of a class of
Schr\"odinger operators, \textit{J.\ Statist.\ Phys.} \textbf{106} (2002), 681--691

\bibitem{lw} Q.-H.\ Liu and Z.-Y.\ Wen, Hausdorff dimension of spectrum of one-dimensional
Schr\"odinger operator with Sturmian potentials, \textit{Potential Anal.} \textbf{20}
(2004), 33--59

\bibitem{loth2} M.\ Lothaire, \textit{Algebraic combinatorics on words}, Cambridge University Press,
Cambridge (2002)

\bibitem{m} H.\ Masur, Interval exchange transformations and measured foliations,
\textit{Ann.\ of Math.} \textbf{115} (1982), 168--200

\bibitem{moo} R.\ Moody (Editor), \textit{The Mathematics of Long-Range Aperiodic Order},
Kluwer, Dordrecht (1997)

\bibitem{mh1} M.\ Morse and G.\ Hedlund, Symbolic dynamics, \textit{Amer.\ J.\ Math.} \textbf{60}
(1938), 815--866

\bibitem{mh2} M.\ Morse and G.\ Hedlund, Symbolic dynamics, II.~Sturmian trajectories, \textit{Amer.\
J.\ Math.} \textbf{62} (1940), 1--42

\bibitem{oprss} S.\ Ostlund, R.\ Pandit, D.\ Rand, H.\ Schellnhuber, and E.\ Siggia,
One-dimensional Schr\"odinger equation with an almost periodic potential, \textit{Phys.\
Rev.\ Lett.} {\bf 50} (1983), 1873--1877

\bibitem{pww} J.\ Peyri\`{e}re, Z.-Y.\ Wen and Z.-X.\ Wen, Polynomes associ\'{e}s aux
endomorphismes de groupes libres, \textit{Enseign.\ Math.} {\bf 39} (1993), 153--175

\bibitem{q} M.\ Queff\'{e}lec, \textit{Substitution Dynamical Systems -- Spectral Analysis},
Springer, Berlin (1987)

\bibitem{r} L.\ Raymond, A constructive gap labelling for the discrete Schr\"odinger operator
on a quasiperiodic chain, Preprint (1997)

\bibitem{rz} R.\ Risley and L.\ Q.\ Zamboni, A generalization of Sturmian sequences:
combinatorial structure and transcendence, \textit{Acta Arith.} {\bf 95} (2000), 167--184

\bibitem{r1} J.\ Roberts, Escaping orbits in trace maps, \textit{Physica A} {\bf 228} (1996),
295--325

\bibitem{rb} J.\ Roberts and M.\ Baake, Trace maps as 3D reversible dynamical systems with
an invariant, \textit{J.\ Stat.\ Phys.} {\bf 74} (1994), 829--888

\bibitem{s1} B.\ Simon, Kotani theory for one dimensional stochastic Jacobi matrices, \textit{Commun.\
Math.\ Phys.} {\bf 89} (1983), 227--234

\bibitem{s2} B.\ Simon, Operators with singular continuous spectrum. I.~General operators, \textit{Ann.\ of
Math.} \textbf{141} (1995), 131--145

\bibitem{simon} B.\ Simon, \textit{Orthogonal Polynomials on the Unit Circle. Part~1. Classical theory},
American Mathematical Society, Providence (2005)

\bibitem{simon2} B.\ Simon, \textit{Orthogonal Polynomials on the Unit Circle. Part~2. Spectral theory},
American Mathematical Society, Providence (2005)

\bibitem{sorspe} E.\ Sorets and T.\ Spencer, Positive Lyapunov exponents for Schr\"odinger operators with
quasi-periodic potentials, \textit{Commun.\ Math.\ Phys.} \textbf{142} (1991), 543--566

\bibitem{s5} A.\ S\"ut\H{o}, The spectrum of a quasiperiodic Schr\"odinger operator, \textit{Commun.\
Math.\ Phys.} {\bf 111} (1987), 409--415

\bibitem{s6} A.\ S\"ut\H{o}, Singular continuous spectrum on a Cantor set of zero Lebesgue
measure for the Fibonacci Hamiltonian, \textit{J.\ Stat.\ Phys.} {\bf 56} (1989),
525--531

\bibitem{s7} A.\ S\"ut\H{o}, Schr\"odinger difference equation with deterministic ergodic
potentials, in \textit{Beyond Quasicrystals} (\textit{Les Houches, 1994}), Springer,
Berlin (1995), 481--549

\bibitem{t} S.\ Tcheremchantsev, Dynamical analysis of Schr\"odinger operators with growing sparse
potentials, \textit{Commun.\ Math.\ Phys.} \textbf{253} (2005), 221--252

\bibitem{van} D.\ Vandeth, Sturmian words and words with a critical exponent,
\textit{Theoret.\ Comput.\ Sci.} \textbf{242} (2000), 283--300

\bibitem{Vee1} W.\ A.\ Veech, Gauss measures for transformations on the space of
interval exchange maps, \textit{Ann.\ of Math.} {\bf 115} (1982), 201--242

\bibitem{Vee2} W.\ A.\ Veech, Boshernitzan's criterion for unique ergodicity of an
interval exchange transformation, \textit{Ergod.\ Th.\ \& Dynam.\ Sys.} {\bf 7} (1987),
149--153

\bibitem{wz} N.\ Wozny and L.\ Q.\ Zamboni, Frequencies of factors in Arnoux-Rauzy
sequences, \textit{Acta Arith.} {\bf 96} (2001), 261--278


\end{thebibliography}
\end{document}